# A Syntactic Adaptive Problem Solver Learning Landscape Structures for Scheduling in Clinical Laboratory

Keyao Wang [1,2] · Bo Liu [2,3*]

[1] *School of Economics and Management, Beihang University, Beijing*

[2] *Academy of Mathematics and Systems Science, Chinese Academy of Sciences, Beijing*

[3] *National Center for Mathematics and Interdisciplinary Sciences, Chinese Academy of Sciences, Beijing*

**Abstract**

This paper attempts to derive a mathematical formulation for real-practice clinical laboratory scheduling, and to present a syntactic adaptive problem solver by leveraging landscape structures. After formulating scheduling of medical tests as a distributed scheduling problem in heterogeneous, flexible job shop environment, we establish a mixed integer programming model to minimize mean test turnaround time. Preliminary landscape analysis sustains that these clinics-orientated scheduling instances are difficult to solve. The search difficulty motivates the search for an adaptive problem solver to reduce repetitive algorithm-tuning work, but with a guaranteed convergence. Yet, under a search strategy, relatedness from exploitation competence to landscape topology is not transparent. Under strategies that impose different-magnitude perturbations, we investigate changes in landscape structure and find that disturbance amplitude, local-global optima connectivity, landscape's ruggedness and plateau size fairly predict strategies' efficacy. Medium-size instances of 100 tasks are easier under smaller-perturbation strategies that lead to smoother landscapes with smaller plateaus. For large-size instances of 200-500 tasks, existing strategies at hand, having either larger or smaller perturbations, face more rugged landscapes with larger plateaus that impede search. Our hypothesis that medium perturbations may generate smoother landscapes with smaller plateaus drives our design of this new strategy and its verification by experiments. Composite neighborhoods managed by meta-Lamarckian learning show beyond average performance, implying reliability when prior knowledge of landscape is unknown.

**Acknowledgements**

This research greatly benefited from multiple field investigations hosted by Associate Professor Mei Wang (Clinical Laboratory at China-Japan Friendship Hospital in Beijing) whose questions and comments were incorporated in the final version. We gratefully acknowledge her fast trust and insight. This research greatly benefited from discussions with Emeritus Professor Michiel Keyzer (SOW-VU, Faculty of Economics and Business Administration, Vrije Universiteit Amsterdam, The Netherlands) whose constructive criticism was incorporated in the final version. During 2007–2009, co-author Bo Liu was a junior researcher at SOW-VU under Prof. Keyzer's supervision and gratefully acknowledges his insight. Bo Liu thanks to and stands with his former SOW-VU Ukrainian colleague, Dr. Viktor Yarovy (Institute for Economics and Forecasting, National Academy of Sciences of Ukraine), with whom Bo Liu benefited greatly from discussions on machine learning and optimization. Special thanks go to Professor Ann Marie Ross (University of Chinese Academy of Sciences) for her hard work on polishing the manuscript. This study was supported in part by Frontier Science Key Research Program, Chinese Academy of Sciences (QYZDB-SSW-SYS020).

---

[*] Corresponding author. Bo Liu, E-mail address: bliu@amss.ac.cn; ORCID: 0000-0002-6615-6601

# 1. Introduction

Clinical laboratories in hospitals must perform medical tests on tube-clinical specimens in a timely way so that medical personnel are enabled to make decisions for patients. Prompt testing of specimens and test results can be urgent. Clinical lab results can support as much as 70% of clinical decisions (Esposito, 2015). In one study, laboratory results were responsible for up to 80% of patient facts required for real-time medicine for critical care in emergency departments (Holland, Smith, & Blick, 2006). Poor laboratory efficiency delays results reporting, significantly impacts follow-up diagnoses, and hazards patients' well-being at life-threatening moments in their lives (Singer, Ardise, Gulla, & Cangro, 2005).

Leveraging clinical laboratory efficiency, especially to decrease test turnaround time, has gone far but without fully implementing optimal scheduling. Innovation in rapid testing methods and in stand-alone high-throughput instruments and their layout has been successfully incorporated, and yet less than optimal scheduling could impede whole laboratory efficacy.

In our field investigations in a clinical laboratory at China-Japan Friendship Hospital in Beijing, we observed how staff laboriously schedule tasks daily with minimal automated assistance. Situations challenged them on how to optimally operate their complex, busy laboratory to achieve higher throughput. Extremely dynamic testing volume, ranging into hundreds of tubes per shift, was common. Round-the-clock lab operation meant multiple shifts needing to be completed. Instruments arrayed in a complex instrument-function topology required staff to figure out a time-saving travel path for each tube passing from pre-analysis to results reporting.

Improving turnaround by means of better clinical laboratory scheduling, and doing so without incurring capital costs from new hardware, seems a worthy startpoint. If sophisticated scheduling tools can automatically arrange a feasible sequence for instrument use with testing tubes, it could improve throughput and prevent medical burnout.

We, consequently, sought a mathematical model and an adaptive problem solver. However, to date, a mathematical formulation expressing the quantifiable relationship among task, resource and capacity in a whole clinical laboratory does not exist (Keskinocak & Savva, 2020). Some steps were confined to scheduling tests on a solo workstation (Aarts, Lindsey, Corkan, & Smith, 1995), or on staff assignments across workstations (Elena, Beraldi, & Conforti, 2006).

As domain-independent problem solving, less dependent on human intervention, is needed for alleviating medical staff's cognitive burden and for guaranteeing convergence, we searched for adaptive problem solvers that could match designated problems with effective algorithms. However, unawareness of relatedness between search strategies and the problem's landscape structures could impede identification of effective strategies (Watson, Barbulescu, Whitley, & Howe, 2002; Watson, Beck, Howe, & Whitley, 2003). A quantifiable characterization of landscape for the clinic-orientated scheduling problem is rather unclear; how landscape topology changes under different strategies, i.e., reflecting instance hardness, remains largely unknown; landscape features that highly influence performance are undeterminable. The central issue was that relatedness from characteristics of the best search strategy to characteristics of landscape is not transparent.

To answer challenges in formulation for a real-practice clinical laboratory scheduling, and to present an adaptive problem solver by leveraging instances' landscape features were the direction of our study. On the way to this destination, we took up three challenging questions.

*Q1: Can an explicit mathematical formulation be derived to quantitatively portray complicated relational structures among process, instrument, and capacity in clinical laboratory?*

*Q2: Can a generic solver be designed under the guidance of explicit landscape structures to make an informed decision on which search strategy performs better on a particular instance?*

*Q3: When an instance is hard for existing strategies at hand, can a new search strategy be defined to be well suited to the instance by leveraging its landscape characteristics?*

Our contributions are four-fold, and Fig. 1 illustrates our approach.

For the first question, we formulated optimal decisions in clinical laboratory in a framework of Distributed scheduling in Heterogeneous and Flexible Job Shop Problem environment (D-HFJSP) (Sect. 3), and established a mixed integer programming model with the objective of minimizing mean of turnaround time (Sect. 4). To our knowledge, this is the first attempt to formulate a whole laboratory scheduling problem and provide a mixed integer programming model.

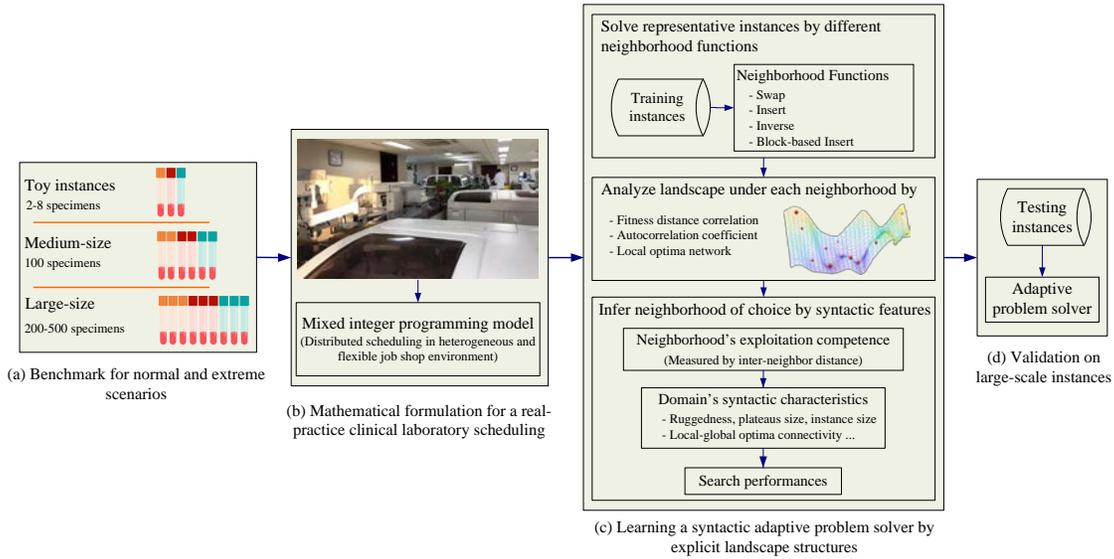

**Fig. 1** A schematic illustration of our approach

For the second question, through offline training, we presented a syntactic adaptive problem solver based on domain's syntactic characteristics (instance size and landscape features) to suggest suitable search strategies for a family of instances (Sect. 5). Here, search strategies were referred to different neighborhood functions in combinatorial space. Three landscape analysis approaches, namely fitness distance correlation (Jones & Forrest, 1995), autocorrelation coefficient (Weinberger, 1990), and local optimal network (Ochoa & Veerapen, 2018) were resorted to in the learning domain's landscape characteristics explicitly. We experimentally discovered that instance hardness under different neighborhood functions can be accurately predicted by instance size, as well as landscape features such as local-global optima connectivity, ruggedness, plateaus size, etc. To our knowledge, our work is the first attempt to understand the nature of landscapes for optimization in healthcare applications.

For the third question, a new neighborhood function for large-size instances was defined by exploiting what we have learned about search spaces. We theoretically proved expected inter-neighbor distances that measured each neighborhood function's competences in exploitation over search space (Appendix A). After experimentally investigated landscape structural changes under each neighborhood function that im-posed perturbations of different magnitudes, we found that, for large-size instances (200 to 500 tasks), the existing neighborhood functions at hand with either larger-perturbations or smaller-perturbations were faced with more rugged landscapes, including multiple scattered funnels (big valleys) with larger attractor or sink plateaus that made the search stagnate. Under moderate perturbations, landscapes were smoother and had smaller attractor plateaus, suggesting large-size instances were easily optimizable under moderate perturbations. This advocated for designing a new neighborhood function, Block-based Insert, to search for optima with moderate momentum and to easily escape from attractor plateaus and traverse multiple funnels (Sect. 6).

Through comprehensive experiments, the syntactic adaptive problem solvers' efficacy was studied (Sect. 7). Our benchmark of 20 toy instances (2, 4, 6, and 8 tasks) and 150 realistic instances (100, 200, 300, 400, and 500 tasks) reflects normal and extreme scenarios. Toy and medium-size instances (100 tasks) were easier under Swap-based strategy that generated smaller perturbations and led to smoother landscape with smaller plateaus. It was validated that Swap achieved the same or better results as Gurobi Optimizer on toy instances but was performing best on medium-size instances. For large-size instances (200 to 500 tasks), the new strategy that allowed a quick step over multiple funnels and plateaus performed best. A generative adaptive problem solver was constructed for comparison, where meta-Lamarckian learning (Ong & Keane, 2004) was relied upon for implicitly learning the landscape structure to suggest the most promising neighborhood function for a single instance at runtime. The proposed syntactic adaptive problem solver surpassed the generative one. On nearly all instances, meta-Lamarckian learning showed beyond average performance, suggesting its reliability when no prior knowledge exists about the problem landscape.

While optimal scheduling is an obvious need in clinical decision support systems, a scheduling model and optimizer that are too complex for physicians is not a workable answer. Our work in applied modeling and generic solver is vital towards automated reasoning and scheduling for machine learning applied to critical clinic testing. The potential for real-world implementation in the clinical laboratory would seem to lie in an explicit mixed integer programming model, together with trained adaptive problem solvers where relatedness from characteristics of the best search strategy to characteristics of landscape is transparent.

## 2. Related Work

To actually derive a mathematical formulation and a generic solver for scheduling in clinical laboratory, our research was extended to a wide scope. Besides investigating scheduling applied to the clinical laboratory, we reviewed topics related to distributed scheduling, adaptive problem solving, and landscape analysis. Herein, we indicate what we adopted from the literature and what we proposed to extend.

**Scheduling in clinical laboratory** research finds test sequences processed on machine with respect to one or more objectives (Davis, 1995). Since Aarts et al. (1995) introduced scheduling theory to clinical chemistry, 25 years of effort seem confined to scheduling a solo workstation, or staff assignments across workstations. Aarts et al. (1995) formulated pipettor dispatch in a chemistry workstation as a no-wait flow shop to minimize turnaround time, and Kim et al. (2006) adopted a job shop model to improve throughput for a portable platform. van Merode et al. (1998) described physician assignments as a mathematical programming model allowing physician preemption. Boyd and Savory (2001) matched by genetic algorithm physicians of different skills with tasks. Elena et al. (2006) presented a multi-period integer programming model to decrease total idle time and the number of physicians needed. To date, mathematical formulation expressing quantifiable relationships among task, resource and capacity in a whole clinical laboratory does not exist (Keskinocak & Savva, 2020). By contrast, we described a mixed integer programming model for this purpose.

**Distributed scheduling** occurs in multiple separate factories or production lines (Behnamian & Ghomi, 2016; Toptal & Sabuncuoglu, 2010). In addition to assignment of jobs to machines and sequencing of jobs on machine, the existence of separate factories has prompted studies on distribution of jobs to factories (Hooker, 2007; J. J. Wang & Wang, 2020). Extensive work modeled each separate factory as a flow shop (Fernandez-Viagas & Framinan, 2015; Hatami, Ruiz, & Andres-Romano, 2013). The work most relevant to our study modeled the separate entity as *homogenous* flexible job shop. Constructive heuristic (Ziaee, 2014), genetic algorithm (Chang & Liu, 2017; De Giovanni & Pezzella, 2010; Lu, Wu, Tan, Peng, & Chen, 2018) and constraint programming (Meng, Zhang, Ren, Zhang, & Lv, 2020) solved the homogenous flexible job shop in distributed scheduling for up to 30 jobs, while an estimation of distribution algorithm (Du, Li, Luo, & Meng, 2021) solved problems with up to 200 jobs. However, the above studies studied benchmark instances rather than real-world cases. To our knowledge, ours is the first report to formulate a *heterogeneous*, flexible job shop in distributed scheduling literature (Lohmer & Lasch, 2020). In healthcare, distributed scheduling has hardly been introduced except for efforts on collaborative operating room scheduling (Roshanaei, Luong, Aleman, & Urbach, 2017).

**Adaptive problem solving,** learns generalized heuristics over a distribution of problems by exploring space of possible heuristics (Gratch & Chien, 1996). It falls within the broader field of automated algorithm selection where by linking behaviour of algorithms to problem features, the (hypothetical) best algorithm from a given set is determined for a problem instance (Rice, 1976). Adaptive problem solving evolved in the path of syntactic, generative, and statistical learning approaches. Informative, interpretable and cheaply computable problem instance features are crucial for learning a good mapping from instances to algorithms (Kerschke, Hoos, Neumann, & Trautmann, 2019). *Syntactic approaches* identify those domain's syntactic structures that influence heuristics' efficacy, such as tightness of constraints (Frost & Dechter, 1994), refinement cost (Kambhampati, Knoblock, & Yang, 1995), or number of feasible serialization orderings (Barrett & Weld, 1994). After finding several efforts had resorted to landscape features as syntactic structures, we elaborated its advance in landscape analysis. *Generative approaches* conjecture heuristics for a single problem based on past problem-solving experiences, examples being PRODIGY (Minton, 1988), case-based reasoning (Rudin et al., 2022; Schmidt, 1998), meta-Lamarckian learning (Ong & Keane, 2004), and adaptive memetic algorithm (Nguyen, Ong, & Lim, 2009). Unlike syntactic and generative approaches, statistical learning approaches (Drake, Kheiri, Özcan, & Burke, 2020; Gratch & Chien, 1996; Gratch & DeJong, 1996; Kool, van Hoof, & Welling, 2019) acknowledge the distribution of problems, and learn policies by imitation (supervised learning) or through experiences (reinforcement learning). Among the challenges that adaptive problem solving encounters is that true distribution of learned and generalized instances cannot be mathematically characterized (Bengio, Lodi, & Prouvost, 2021). We did not address these open questions and, instead, charted an alternative route by proposing an adaptive problem solver embodying the syntactic approach. Our rationale was that by leveraging the landscape structures explicitly, the most suitable strategy would make the search easier, and consequently can be learned for an instance.

**Landscape analysis**, since it reflects topology of combinatorial search space (Malan & Engelbrecht, 2013; Reeves, 1999; Reidys & Stadler, 2002), helps identify landscape features that highly influence search performances, that characterize problem hardness under different strategies, and that motivate new strategy construction (Daolio, Liefooghe, Verel, Aguirre, & Tanaka, 2017). Hutter et al. (2014) modeled algorithm runtime as a function of problem-specific features, where autocorrelation coefficient quantified TSP's land-scape ruggedness. Watson et al. (2005) found that mean distance between random local optima and the nearest optima highly correlates with problem hardness for tabu search on job shops. Streeter and Smith (2006) characterized job shop's landscape as a big valley for low ratio of job to machine, and as many separate big valleys for high ratio of job to machine. Franzin and Stützle (2023) answered the open question whether "to cool or not" in simulated annealing by observing whether a certain landscape property changed or held constant throughout the search space. Several studies showed how hypothesis on landscapes could motivate design of new algorithms. Merz and Katayama (2004) upon seeing that local optima for quadratic programming problems were concentrated in a small fraction of the solution domain, designed large-jump mutation to traverse the basins of attraction. Zhang (2004) estimated and utilized backbone structures to guide search for Boolean satisfiability. Barbulescu et al. (2006) found plateaus dominated search spaces for satellite control network scheduling, and constructed algorithms making larger perturbations to solutions. However, characterization of landscapes for distributed scheduling is rather unclear to-date; how landscape structure changes under different neighborhood functions remains unknown. Once we analyzed landscapes for a family of clinical laboratory scheduling instances, comprehensively learned that domain's syntactic structures and acknowledged relatedness between landscape structure and search strategies, we could then infer the choice of a neighborhood function for the adaptive problem solver, but also define a new neighborhood function.

An explicit mathematical formulation, portraying the complicated relations among task, resource and capacity in a whole clinical laboratory, does not exist. An adaptive problem solver taking advantage of explicit landscape features as domain's syntactic structures is highly needed. These motivate us to present a sophisticated, automated scheduling tool, thereby achieving the best clinical practice and enhancing healthcare quality in a time-critical situation.

## 3. Workflows and Modeling Aspects for Scheduling in Clinical Laboratory

We introduce workflows for processing medical specimens in the clinical laboratory and characterize the scheduling model.

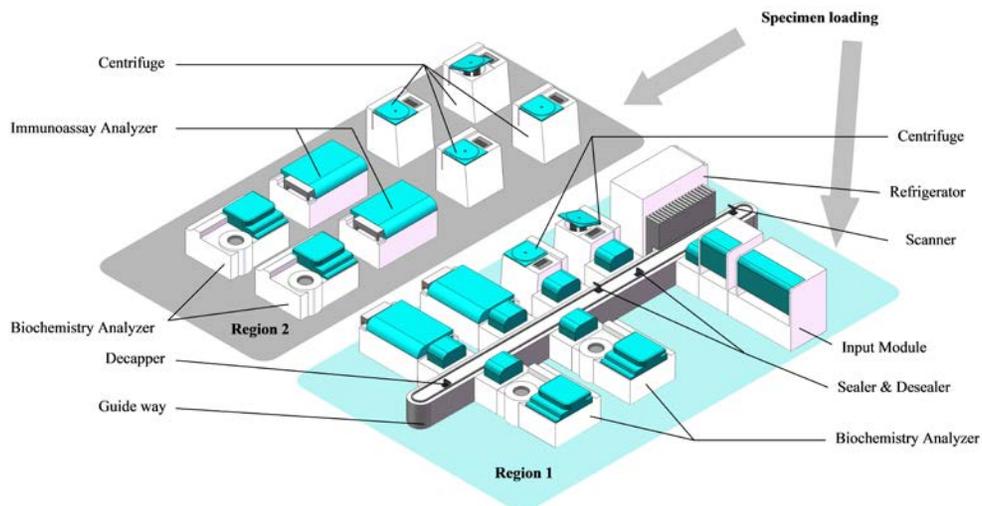

**Fig. 2** Representation of workflows for processing medical specimens, and instrument layout in a real clinical laboratory

## 3.1 Workflows for Processing Medical Specimens in Clinical Laboratory

A clinical laboratory offers a variety of testing services. We chose two exceedingly common ones, blood biochemical and immunologic tests. In Fig. 2, a 4-stage workflow for processing blood specimens, involving testing order recognition, specimen preparation, analysis, and results validation and reporting, is illustrated in a schematic of a real laboratory's workflow and instrument layout. Upon receiving specimens from outpatient/inpatient drawing centers, lab staffs distribute specimens to one of two lines, Region 1 or Region 2. Automated scanners (Region 1) or staff members (Region 2) recognize test requests by reading tube barcodes, and group them into lots. Next, blood specimens will be prepared by centrifuges to separate serum and plasma, and by decapper or staff to remove lids. Third, tubes grouped in batches are loaded to biochemistry or immunoassay analyzers. Last, testing results from analyzers are validated and reported by physicians. Afterwards, tubes are sealed and stored in refrigerator for later use.

We interviewed the staff at Clinical Laboratory of China-Japan Friendship Hospital in Beijing to obtain data needed for scheduling. Lower and upper bounds of processing time, maximum capacity, and number of parallel instruments are listed in Table 1. Maximum capacity refers to the upper bound of the instrument that can simultaneously process tubes.

**Table 1** Operations, processing time, capacity and quantity of instruments/physicians in a real-practice clinical laboratory

| Region | Instrument / Staff | Operation | Processing time (second) | Maximum Capacity | Number of parallel instruments |
|---|---|---|---|---|---|
| 1 | Centrifuge | Centrifugation | [480, 600] | 84 | 2 |
| 1 | Decapper | Decapping | 2 | 1 | 1 |
| 1 | Biochemistry Analyzer | Biochemical test | [480, 600] | 84 | 2 |
| 1 | Immunoassay Analyzer | Immunologic test | [1080, 1800] | 84 | 2 |
| 1 | Staff | Result validation & reporting | [4, 6] | 1 | 1 |
| 2 | Centrifuge | Centrifugation | [300, 360] | 32 | 4 |
| 2 | Staff | Decapping | [4, 6] | 1 | 1 |
| 2 | Biochemistry Analyzer | Biochemical test | [300, 720] | 60 | 2 |
| 2 | Immunoassay Analyzer | Immunologic test | [900, 2700] | 48 | 2 |
| 2 | Staff | Result validation & reporting | [4, 6] | 1 | 1 |

## 3.2 Modeling Aspects for Scheduling in Clinical Laboratory

A variety of aspects can be considered when developing scheduling models. The roadmap proposed by Mendez et al. (2006) was used; it considers process topology, equipment issues, and time and demand-related constraints. See Mendez et al. (2006) for details. The following modeling aspects were considered.

*Process topology*: Modelled as Distributed scheduling in Heterogeneous and Flexible Job Shop Problem environment, it is considered to be *distributed* in that there are two separate lines, *heterogeneous* in that lines can differ in instruments' quantity, capacity and speed (processing time), and *flexible job shop* in that each line has parallel instruments at certain stages and each tube has its own pre-determined operational sequence.

*Equipment connectivity and assignment*: Each tube can only be assigned to one line. Once assigned, tubes cannot transfer across lines, which imposes hard constraints on equipment allocation.

*Specimen transportation*: Movement of specimens between instruments is performed by guide way in Region 1 or by human action in Region 2. As transport time between instruments is rather short as compared with specimens' processing time, transport time is neglected in this study.

*Batch size*: Batch size cannot exceed an instrument's maximum capacity.

*Processing time*: Processing time is considered to be fixed. Processing time for a tube in a batch equals processing time for the batch.

*Demand patterns*: Requested testing orders in short-term scheduling are assumed to be fixed.

*Changeovers*: Consumables (including reagents and probes) are assumed to be sufficient for short-term scheduling. Thus, changeover time is not considered.

*Resource constraints*: Resource constraints focus on instruments while other constraints (including labor, power, reagents, and other utilities) are not considered.

*Time constraints*: Time constraints such as non-working periods on the weekend or maintenance periods are not considered in short-term scheduling.

*Costs*: In our field investigation, because the mean in-lab to reporting turnaround time (MTAT) is crucial, minimization of MTAT is considered.

The overall situation is characterized as a distributed scheduling problem in which each separate line is a flexible job shop. They are heterogeneous. Considering batch machines, the objective is to minimize mean test turnaround time. The problem is abbreviated as D-HFJSP. According to the $\alpha/\beta/\gamma$ notation for classifying scheduling problems (Pinedo, 2012), D-HFJSP is denoted as $D_F, FJ_c/batch/MTAT$. It is NP-hard since it is more complex than the NP-hard flexible job shop (Pezzella, Morganti, & Ciaschetti, 2008).

## 4. Problem Formulation

In this section, we establish a mixed integer programming model for D-HFJSP. The indices, parameters and variables are explained in Table 2.

**Table 2** Indices, parameters and variables in mixed integer programming model

| **Indices** | |
|---|---|
| $i$ | Index of specimens |
| $j$ | Index of operations |
| $l$ | Index of job shop lines |
| $k$ | Index of machines in job shop line |
| $d$ | Index of batches |
| $r$ | Position of batch in the processing sequence on machine |
| **Parameters** | |
| $n$ | Number of specimens |
| $J_i$ | The $i$-th specimen |
| $o_i$ | Number of operations of specimen $J_i$ |
| $O_{i,j}$ | The $j$-th operation of specimen $J_i$ |
| $f$ | Number of job shop lines |
| $F_l$ | The $l$-th job shop line |

| $m_l$ | Number of machines in job shop line $F_l$ |
|---|---|
| $M_{l,k}$ | The $k$-th machine in job shop line $F_l$ |
| $c_{l,k}$ | Capacity of machine $M_{l,k}$ in job shop line $F_l$ |
| $W_{i,j}^{l,k}$ | Binary parameter that takes value 1 if machine $M_{l,k}$ is capable of executing operation $O_{i,j}$ of specimen $J_i$; 0 otherwise |
| $P_{i,j}^{l,k}$ | Processing time of specimen $J_i$ for operation $O_{i,j}$ on machine $M_{l,k}$ in job shop line $F_l$ |
| $\Gamma$ | A sufficiently large positive number |
| **State Variables** | |
| $B_d^{l,k}$ | The $d$-th batch processed on machine $M_{l,k}$ in job shop line $F_l$ |
| $T_d^{l,k}$ | Processing time of batch $B_d^{l,k}$ |
| $E_{i,j}$ | Available time of operation $O_{i,j}$ of specimen $J_i$ |
| $C_r^{l,k}$ | Completion time of the $r$-th batch on machine $M_{l,k}$ |
| $TAT_i$ | Turnaround time of specimen $J_i$ |
| **Decision Variables** | |
| $X_i^l$ | Binary variable, taking value 1 if specimen $J_i$ is processed in job shop line $F_l$; 0 otherwise |
| $Y_{i,j,d}^{l,k}$ | Binary variable, taking value 1 if the $j$-th operation $O_{i,j}$ of specimen $J_i$ belongs to the $d$-th batch $B_d^{l,k}$ processed on the $k$-th machine $M_{l,k}$; 0 otherwise |
| $Z_d^{l,k,r}$ | Binary variable, taking value 1 if the $d$-th batch $B_d^{l,k}$ is processed in the $r$-th position on machine $M_{l,k}$; 0 otherwise |

The model is represented as

$$\min\left\{MTAT = \frac{1}{n} \cdot \sum_{i=1}^{n} TAT_i\right\} \tag{1}$$

Subject to:

$$\sum_{l=1}^{f} X_i^l = 1, \ \forall i \tag{2}$$

$$\sum_{k=1}^{m_l} \sum_{d=1}^{n} W_{i,j}^{l,k} \cdot Y_{i,j,d}^{l,k} = X_i^l, \ \forall i,j,l \tag{3}$$

$$\sum_{l=1}^{f} \sum_{k=1}^{m_l} \sum_{d=1}^{n} Y_{i,j,d}^{l,k} = 1, \ \forall i,j \tag{4}$$

$$\sum_{i=1}^{n} \sum_{j=1}^{o_i} Y_{i,j,d}^{l,k} \le c_{l,k}, \ \forall d,k,l \tag{5}$$

$$\sum_{d=1}^{n} Z_d^{l,k,r} \le 1, \ \forall k,l,r \tag{6}$$

$$\sum_{r=1}^{n} Z_d^{l,k,r} = 1, \ \forall d,k,l \tag{7}$$

$$T_d^{l,k} = \max_{i,j}\left\{Y_{i,j,d}^{l,k} \cdot P_{i,j}^{l,k}\right\}, \ \forall i,j,d,k,l \tag{8}$$

$$C_1^{l,k} \ge T_d^{l,k} + \Gamma \cdot \left(Z_d^{l,k,1} - 1\right), \ \forall d,k,l \tag{9}$$

$$C_r^{l,k} \ge C_{r-1}^{l,k} + T_d^{l,k} + \Gamma \cdot \left(Z_d^{l,k,r} - 1\right), \ \forall d,k,l, r=2,\ldots,n \tag{10}$$

$$C_r^{l,k} \ge E_{i,j} + T_d^{l,k} + \Gamma \cdot \left(Y_{i,j,d}^{l,k} + Z_d^{l,k,r} - 2\right), \ \forall d,i, j=1,\ldots o_i, k,l,r \tag{11}$$

$$E_{i,(j+1)} \ge C_r^{l,k} + \Gamma \cdot \left(Y_{i,j,d}^{l,k} + Z_d^{l,k,r} - 2\right), \ \forall d,i, j=1,\ldots,(o_i-1), l,k,r \tag{12}$$

$$E_{i,1} \ge 0, \ \forall i \tag{13}$$

$$C_r^{l,k} \ge 0, \ \forall k,l,r \tag{14}$$

$$TAT_i \ge C_r^{l,k} + \Gamma \cdot \left(Y_{i,o_i,d}^{l,k} + Z_d^{l,k,r} - 2\right), \ \forall d,i,l,k,r \tag{15}$$

$$X_i^l \in \{0, 1\}, \ \forall i, l \tag{16}$$

$$Y_{i,j,d}^{l,k} \in \{0, 1\}, \ \forall d, i, j, k, l \tag{17}$$

$$Z_d^{l,k,r} \in \{0, 1\}, \ \forall d, k, l, r \tag{18}$$

- Objective (1) is to minimize mean test turnaround time for all specimens $J_i \in \{J_1, ..., J_n\}$.

- Constraint (2) guarantees that each specimen $J_i$ can only be assigned to one job shop line $F_l$, $F_l \in \{F_1, ..., F_f\}$. Once distributed, specimens are not allowed to transfer across lines.

- Constraints (3) and (4) guarantee that the $j$-th operation $O_{i,j}$, $O_{i,j} \in \{O_{i,1}, ..., O_{i,o_i}\}$ for specimen $J_i$ can only be allotted to one of the batches on machine $M_{l,k}$ in the dedicated line, and they guarantee that $O_{i,j}$ is assigned to the eligible machine $M_{l,k}$, $M_{l,k} \in \{M_{l,1}, ..., M_{l,m_l}\}$ where $W_{i,j}^{l,k} = 1$.

- Constraint (5) ensures that lot size does not exceed the machine's capacity $c_{l,k}$.

- Constraints (6) and (7) determine each lot's position in the machine's processing sequence by guaranteeing each lot can only be placed in one position, and each position can hold at most one lot.

- Constraint (8) implies that the processing time of the $d$-th batch is equal to the longest processing time of the specimen in the batch.

- Constraints (9) and (10) imply that the completion time of the $r$-th batch on machine $M_{l,k}$ should be not less than the completion time of the $(r-1)$-th batch on machine $M_{l,k}$ plus the processing time of the $r$-th batch.

- Constraint (11) specifies that batch completion time is not less than the available time of all specimens allotted to the batch plus the processing time of the batch, i.e., the batch's processing cannot start on the current machine until all specimens of that batch complete their preceding operations and are released.

- Constraint (12) implies that no specimen becomes available for processing on the machine of the next processing stage until the batch to which it belongs has been processed.

- Constraint (13) implies that at time 0 all specimens are available for processing.

- Constraint (14) sets batch completion times to be non-negative.

- Constraint (15) defines the turnaround time for specimen $J_i$.

To conclude, we derived an explicit mixed integer programming model to portray the complicated relationship in clinical laboratory. By solving it, optimal decisions are obtained that boil down to four consecutive decisions: distribution of tubes to lines, assignment of tubes to instruments, sequencing, and time of tubes on instrument. The next section illustrates the solver's main features.

## 5. Syntactic Adaptive Problem Solver

In this section, we present the syntactic adaptive problem solver, and explain its components related to landscape analyzers, neighborhood functions, solution representation and evaluation, and three commonly used combinatorial optimization algorithms.

## 5.1 Framework and Procedure for the Syntactic Adaptive Problem Solver

It is generally accepted that the landscape depends not only on the objective function and instance data, but also on the algorithm components used, in particular the neighborhood function (Franzin & Stützle, 2023). A search for a generic problem solver, to match designated instances with effective neighborhood function, is non-trivial and challenging mainly because landscape topology that reflects instance hardness under a neighborhood function is unknown. Therefore, in our adaptive problem solver, we focus on discovering the impact of instance data, neighborhood functions on landscape features. In this way, the revealed transparent influence relationships are expected to help match neighborhood functions for different instances, leading to better search performance. Since only one objective function is considered in this study, it is not necessary to consider the effect of the objective function on the landscape features; if there are multiple objective functions, it is required to examine their impact on landscape features.

A schematic of the adaptive problem solver based on domain's syntactic characteristics is illustrated in Fig. 3. The adaptive problem solver operated as follows.

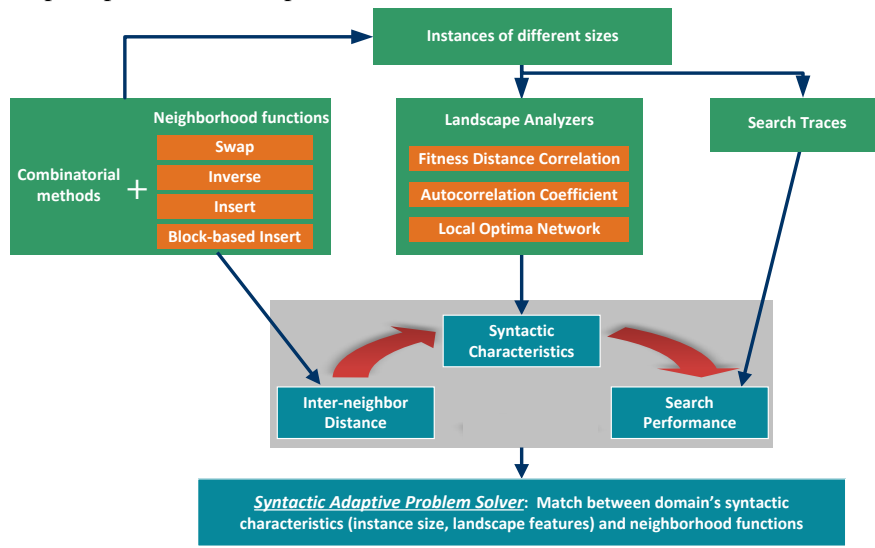

**Fig. 3** A schematic of the adaptive problem solver based on domain's syntactic characteristics

First, we theoretically proved the expected inter-neighbor distance that measures the competence in exploitation over the solution space for each neighborhood function (Sect. 6.2 and Appendix A). The inter-solution distance is measured by Job Precedence Rule-based distance, JPR distance (Liu et al., 2013). It counts the number of job pairs with identical elements but different precedence in two solutions. JPR distance metric is detailed in Appendix D.

Second, representative instances are solved by the same combinatorial optimization algorithm but with different neighborhood functions, and search traces are recorded for each representative instance under each neighborhood function. We have following notes. (*i*) In our preliminary investigation, the running time increased exponentially with instance size. We partitioned instances into subgroups by size (100, 200, 300, 400, and 500), simply selecting instances from each subgroup as representative. (*ii*) Three combinatorial optimization algorithms were used for experiments. We selected two representative algorithms, Simulated Annealing (SA) and Scatter Search (SS), from single-solution-based and population-based metaheuristics, respectively. We also used a SA variant, the Fixed Temperature algorithm (FTA) where the temperature remains constant throughout the search. Our main goal is to observe the effect of different neighborhood functions on the landscape structures, and consequently on the search performances. In addition to the aforementioned three commonly used combinatorial optimization methods, other methods can also be used as the carrier of the neighborhood functions. (*iii*) We solved representative instances of each subgroup with each combination of algorithms and neighborhood functions.

Third, through three landscape analyzers, namely fitness distance correlation (Jones & Forrest, 1995), autocorrelation coefficient (Weinberger, 1990), and local optimal network (Ochoa & Veerapen, 2018) (Sect. 5.2), instances' landscape characteristics were explicitly calculated, revealing how landscapes change under different neighborhood functions that imposed perturbations of different magnitudes (Sect. 6).

Fourth, we examined the chain relatedness starting from search strategy's exploitation competence over search space (measured by inter-neighbor distance), going through landscape topology, and then to search performance. Through experimental observations, we determine landscape characteristics that highly influence search performance. We relate characteristics of the best search strategy to characteristics of landscape, and infer the choice of a strategy that made search easier for each subgroup. We have following note. (*i*) We revealed that solution quality or instance hardness under different neighborhood functions could be accurately predicted by instance size, as well as landscape features such as local-global optima connectivity, ruggedness, plateaus size, etc. (*ii*) The main challenge in designing a syntactic adaptive problem solver—which landscape features to choose as syntactic features to accurately reflect instance hardness under neighborhood functions—is addressed.

Fifth, through the above offline training, we set up an adaptive problem solver based on domain's syntactic characteristics (instance size and local-global optima connectivity, ruggedness, plateaus size) to infer the most suitable neighborhood function. Once a syntactic adaptive problem solver is trained, when solving new instances, the most appropriate neighborhood function is determined for a combinatorial algorithm based on the learned transparent relationship between domain's syntactic characteristics and neighborhood functions.

Finally, the syntactic adaptive problem solver was tested on large-scale instances, and statistical analysis was performed to validate its performances (Sect. 7).

Next, we explain main components used in the adaptive problem solver.

**5.2 Landscape Analyzers**

Characteristics of landscape are critical indicators of search space topology, determining instance hardness under different neighborhood functions. We use three analyzers, namely fitness distance correlation, autocorrelation coefficient, and local optima network. Through these analyzers, landscape can be quantitatively characterized in terms of connectivity from local to global optima, ruggedness, and the size of plateaus. These landscape characteristics are recognized as the domain's syntactic characteristics, and we hypothesize they influence search efficacy.

5.2.1 Fitness Distance Correlation

Fitness distance correlation (Jones & Forrest, 1995) reveals connectivity from local optima to global optima, reflecting search difficulties towards global optima. Large fitness distance correlation indicates strong local-global optima connectivity, implying easier movement from local optima to global optima (or the best-so-far solutions), as estimated by their correlation with respect to their distance, as in

$$FDC = \frac{\sum_{i=1}^{m}(f(x_i) - \overline{f})(d(x_i) - \overline{d})}{m\sigma_f \sigma_d} \quad (19)$$

where $m$ is the size of set of solutions $x_i$, $f(x_i)$ is the objective value, $d(x_i)$ is the distance from $x_i$ to the global optima (or the best-so-far solution) and is measured by JPR distance metric (see Appendix D), $\overline{f}$ and $\overline{d}$ are the mean of $f(x_i)$, and $d(x_i)$, $\sigma_f$ and $\sigma_d$ are the standard deviations of $f(x_i)$ and $d(x_i)$.

5.2.2 Autocorrelation Coefficient

Autocorrelation coefficient (Weinberger, 1990) reveals landscape ruggedness, reflecting search difficulties around its neighbors. A zero autocorrelation value sustains a rugged landscape. Conversely, a one autocorrelation value suggests a smooth landscape. It is estimated by the correlation of solutions that are several steps apart along a random walk.

$$AC(s) = \frac{\sum_{t=1}^{m-s}(f(x_t)-\overline{f})(f(x_{t+s})-\overline{f})}{\sigma_f^2(m-s)} \qquad (20)$$

where $m$ is random walk length, $x_t$ is the $t$-th solution along the random walk, and $s$ is the number of steps apart from two sequences along the random walk.

5.2.3 Local Optima Network

Local optima network is modeled with nodes (local optima) connected by directed, weighted edges (transition probability from a node to another) (Ochoa & Veerapen, 2018). The network reflects risks of getting stuck in local optima and difficulties in traversing basins of attraction. Via local optima network, Daolio et al. (2014) characterized flow shop's landscape in terms of network size, clustering coefficient, transition probabilities, link heterogeneity, path length, and mixing patterns. These features fairly explained instance hardness under iterated local search and predicted performances.

To our knowledge, for modeling real-practice scheduling problems, local optima network has not yet been used. In view of its successful application in flow shop, we believe it has potential in speculating the landscapes although our clinics-orientated scheduling problems are more complex. Table 12 in Appendix E gives a brief explanation about concepts in local optima network, and Fig. 4 illustrates the concepts. To visualize local optima network, we develop an algorithm to obtain nodes and edges (see Appendix F), and then adopt a force-directed method (Ochoa & Veerapen, 2018) to assign nodes to positions in a metric space.

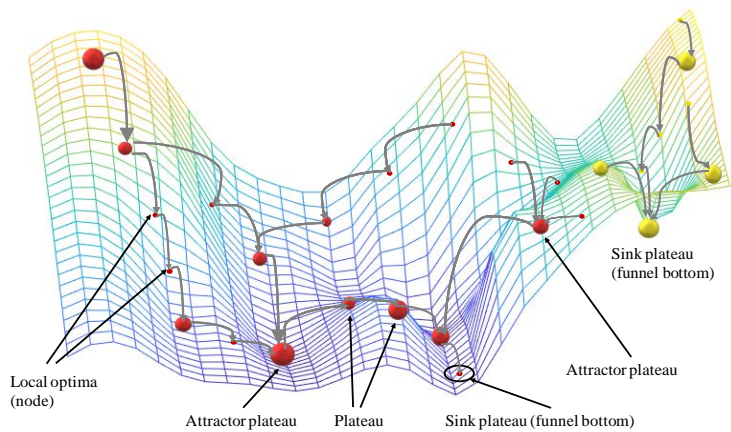

**Fig. 4** An illustrative local optima network on a landscape for a minimization problem, showing the two funnels, one connected to red colored nodes and the other to yellow-colored nodes

## 5.3 Neighborhood Functions

In the adaptive problem solver, a total of four neighborhood functions for combinatorial space are used, namely Insert, Swap, Inverse, and Block-based Insert.

1) Insert (INS): Choose randomly two distinct elements from vector of specimen sequence which is a $n$-length permutation. Insert the back one before the front. As illustrated in Fig. 5a, specimens 1 and 5 are selected, and specimen 5 (the latter) is inserted in front of specimen 1.

2) Swap (SWP): Select randomly two distinct elements from the sequence and swap. Specimens 1 and 5 are selected and swapped, as shown in Fig. 5b.

3) Inverse (INV): Invert the subsequence between two different random positions in the sequence. Specimens 1 and 5 are selected. The inverse, the partial permutation between them, is shown in Fig. 5c.

4) Block-based Insert (INB): Choose randomly two distinct blocks that are generated in Sect. 5.6.1. Insert the back block before the front one without changing specimens' priority within each block. As depicted in Fig. 5d, Blocks 2 and 3 are selected, and Block 3 is inserted in front of Block 2 while keeping specimens' priority within each block unchanged. Another illustration is shown in Fig. 5e where Block 3 is in front of Block 2, since the first element of Block 3 is before the first element of Block 2. All elements of Block 2 are inserted before the first element of Block 3.

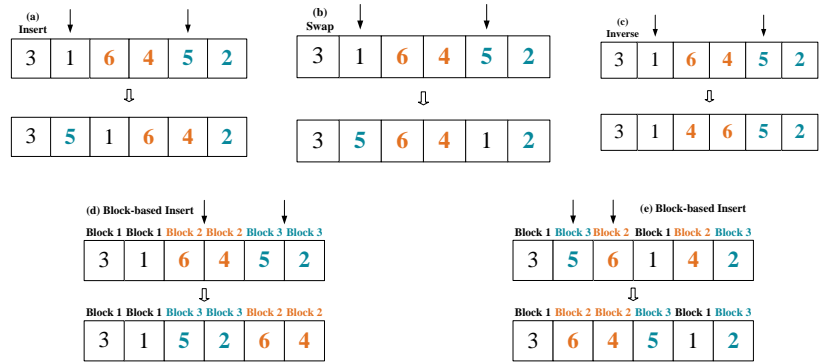

**Fig. 5** Illustration of Neighborhood functions: (a) Insert, (b) Swap, (c) Inverse, (d, e) Block-based Insert

Our initiative to design Block-based Insert is to impose medium perturbations to generate a smooth landscape for large-size instances (200-500 tasks). We find in Sect. 6.3-6.5 that for large-size instances, existing neighborhood functions (Insert, Swap, Inverse) at hand generate either larger or smaller perturbations, resulting in rugged landscapes with large plateaus, thereby impeding search. Our hypothesis is that medium perturbations may generate smoother landscapes with smaller plateaus. Compared with Insert, Swap and Inverse, the Block-based Insert has a moderate inter-neighbor distance; see Sect. 6.2. Under Block-based Insert, landscapes are smoother and have smaller attractor plateaus, suggesting large-size instances are easily optimizable. We will experimentally investigate changes in landscape structure under each neighborhood function that imposed perturbations of different magnitudes, and show how search performances are affected in Sect. 7.

## 5.4 Solution Representation and Evaluation

We provide the solution formulation and evaluation—a common component—required in the optimization algorithms of this study before introducing them.

### 5.4.1 Solution Representation

We adopt the most common permutation-based encoding scheme (Liu, Wang, & Jin, 2007) to represent solution in the form of vector of specimen sequence (VSS). The element's position represents processing precedence among all elements.

### 5.4.2 Solution Evaluation

To transform the VSS-based solutions into feasible schedules (sequencing and timing of operations to machines), a novel and efficient dispatching heuristic named First Available Batch Machine (FABM) is proposed to (1) distribute specimens to lines, (2) group specimens into batches, and (3) assign batches to eligible machines. FABM is illustrated in Fig. 6. Note that the multi-list partial permutation widely used in the distributed permutation flowshop (Naderi & Ruiz, 2014) cannot handle the situations encountered in this study, such as batch machine, parallel machine, and different pre-determined operational sequence for each job.

**BEGIN**
 Input a vector of specimen sequence (VSS).
 (1): Distribute specimens to lines. In order left to right, select one element at a time from VSS, and implement Step (1).
  (1.1): For selected specimen, determine a set of eligible machines that can process its first operation.
  (1.2): Determine a set of available machines from the eligible set in Step (1.1).
   (1.2.1): If the set of available machines is not empty, go to Step (1.3).
   (1.2.2): If the set of available machines is empty, suggesting all eligible machines are occupied, then according to Constraints (9) – (11), an eligible machine that has an earliest available time can be determined. Put the specimen on that machine's waiting list.
  (1.3): Determine the machine from the available set obtained in Step (1.2).
   (1.3.1): If there is only one machine in an available set, assign the specimen to this machine.
   (1.3.2): If there are multiple machines, assign the specimen to the first machine whose capacity is closest to its maximum capacity after adding the specimen into its processing batch.
   (1.3.3): If there are multiple machines with the same remaining capacity, randomly assign specimen to any machine.
  (1.4): Specimen is distributed to the job shop line that contains the machine.
  *Remark*: Hereafter, the specimen's remaining operations can only be processed on machines in the same job shop line, and it is not allowed to transfer across lines. This complies with Constraint (2).
 (2): Group specimens for batch processing on machines.
  (2.1): For each machine, specimens assigned to this machine are sorted in the ascending of arriving time that is calculated from Constraint (12). If more than one specimen has the same arriving time, then they are ranked according to the precedence specified in VSS.
  (2.2): Split ranked specimens into batch(es) by size equal to machine capacity. If the number of remaining specimens is less than the machine capacity, these specimens are regarded as a single batch.
  (2.3): For each batch assigned to the machine, the machine can only begin processing once this batch's specimens are all available. Their arriving times are computed by Constraint (11).
  *Remark:* FABM is proposed for batch processing machines (e.g., centrifuge, biochemistry analyzer and immunoassay analyzer), but is also applicable for non-batch machines simply by limiting machine capacity to one (e.g., decapper).
 (3): Dispatch the successive operation for each specimen on machine according to Steps (1.1)–(1.3), and group specimens for batch processing on that machine according to Step (2).
 (4): Repeat Step (3) until all operations of all specimens are assigned.
 (5): Output the sequencing and timing of operations to machines, that is, the schedule.
**END**

**Fig. 6** First Available Batch Machine (FABM) generates feasible schedules

### 5.4.3 Instance to Explain the First Available Batch Machine

We provide an instance to explain FABM. In the illustrative case, 6 specimens are processed by two job shop lines. Specimens 1-3 are for biochemical tests and specimens 4-6 for immunologic tests. Processing time for each operation by eligible machines is listed in Table 10 in Appendix B. Time is sampled within the lower and upper bounds of processing time in Table 1. The sign "-" means the machine is ineligible for the operation. For instance, a value of 545 at the intersection of row $O_{1,3}$ and column $(M_{1,4}, M_{1,5})$ means that the 4th and 5th machines are both eligible for the 3rd operation of specimen $J_1$ and the processing time is 545 seconds. For illustrative purpose, we set maximum capacity to be [2, 1, 2, 2, 1] for centrifuges, decappers, biochemistry analyzers, immunoassay analyzers, and results validation & reporting, respectively.

Assuming VSS to be [3, 1, 6, 4, 5, 2], for the first specimen $J_3$, its first operation $O_{3,1}$ could be randomly assigned to any machine in eligible machine set $\{M_{1,1}, M_{1,2}, M_{2,1}, M_{2,2}\}$ according to Steps (1.1)-(1.3) in FABM. More specifically, according to Step (1.3.3), $O_{3,1}$ can be assigned to any machine because the four machines' remaining capacities are the same. Here, $O_{3,1}$ is assigned to $M_{2,1}$. For the second specimen $J_1$, its first operation $O_{1,1}$ is assigned to $M_{2,1}$ according to Steps (1.1), (1.2.1) and (1.3.2), so that $M_{2,1}$ can reach its maximum capacity and begin processing $\{O_{3,1}, O_{1,1}\}$ as soon as possible. Til now, specimens $J_3$ and $J_1$ have been assigned to a machine in line $F_2$, and their remaining operations can only be processed on machines in line $F_2$ according to Step (1.4). For the remaining specimens in VSS, by implementing Step (1) on an element-by-element basis, $O_{6,1}$ and $O_{4,1}$ are assigned to $M_{1,2}$, and $O_{5,1}$ and $O_{2,1}$ are to $M_{2,2}$. To sum up, specimens $\{J_6, J_4\}$ and $\{J_3, J_1, J_5, J_2\}$ are distributed to $F_1$ and $F_2$, respectively.

According to Step (2) in FABM, $J_3$ and $J_1$ are grouped as the first batch on machine $M_{2,1}$ for batching processing, and the batch's processing time is the maximum processing time of $O_{3,1}$ and $O_{1,1}$ on $M_{2,1}$, i.e., 332 seconds. Following Step (2), $J_6$ and $J_4$ are grouped as the first batch on $M_{1,2}$, and batch processing time is 564 seconds. $J_5$ and $J_2$ are grouped as the 1st batch on $M_{2,2}$, and the processing time is 356 seconds.

The successive operation for each specimen is dispatched to machines by Step (3). For the second operation of each specimen, only one eligible machine in each job shop line, i.e., $M_{1,3}$ and $M_{2,3}$, is available. Once $J_3$ and $J_1$ are released from batching machine $M_{2,1}$, they arrive at the successive eligible machine $M_{2,3}$ simultaneously. Since $M_{2,3}$'s maximum capacity is 1, $O_{3,2}$ is first processed according to Steps (2.1)-(2.2). $J_1$ is on the waiting list of $M_{2,3}$, and $O_{1,2}$ can be processed after $O_{3,2}$ is finished. By the same means, $\{O_{5,2}, O_{2,2}\}$ and $\{O_{6,2}, O_{4,2}\}$ are dispatched and processed by machine $M_{2,3}$ and $M_{1,3}$, respectively. We repeat Step (3) until all operations for all specimens are assigned. By implementing FABM, we obtain the schedule, i.e., the sequencing and timing of the machines' operations. Fig. 7 illustrates the resultant Gantt chart of the schedule. Decision variables ($X_i^l$, $Y_{i,j,d}^{l,k}$, and $Z_d^{l,k,r}$) corresponding to this schedule are valued in Table 11 in Appendix C.

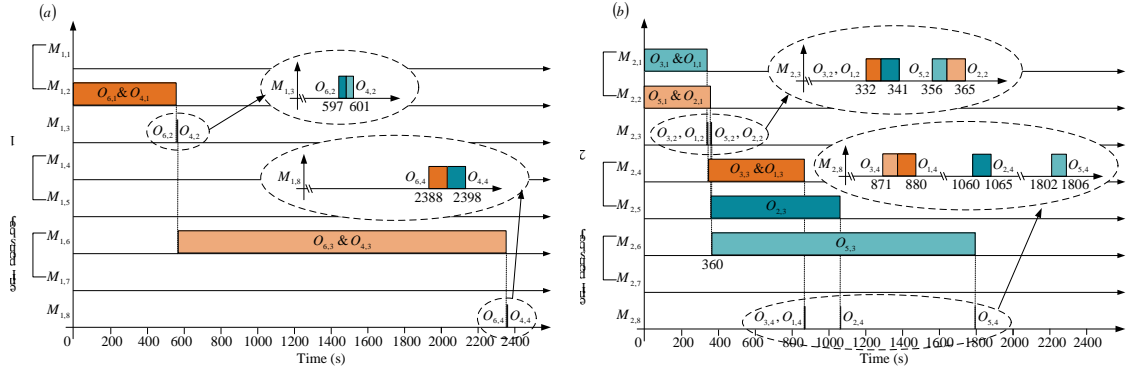

**Fig. 7** Gantt chart depicting the assignment of operations to machines, sequencing of operations, and timing of operations

### 5.5 Simulated Annealing (SA) and Fixed Temperature Algorithm (FTA)

In combinatorial optimization, simulated annealing usually achieves state-of-the-art performance despite its simplicity (Kirkpatrick, Gelatt, & Vecchi, 1983). The simulated annealing starts from an initial state. In our case, the initial state is a random sequence of all pending specimens. A new candidate is randomly sampled from its neighbors, causing a change in objective value. The new candidate is accepted to replace the current state with an acceptance probability related to temperature and change in objective value. At each temperature, the simulated annealing deploys the neighborhood to perform perturbations with the size of Metropolis sampling (Metropolis, Rosenbluth, Rosenbluth, Teller, & Teller, 1953). Once the sampling size is reached, the temperature is updated according to annealing schedule. The above search process is repeated until the temperature drops to the specified minimum temperature. Since the probability to accept an inferior candidate gradually decreases with the temperature drops, simulated annealing has potential to escape local optima, and has a speed-up convergence.

Parameters, such as initial temperature, annealing schedule, acceptance probability, size of Metropolis sampling, and neighborhood function are related to its search efficacy (Ball, Branke, & Meisel, 2018). We set an initial temperature $T_0$ by trial and error, use an exponential annealing schedule $T_k = \lambda \cdot T_{k-1}$ where $\lambda \in (0, 1)$ is the cooling coefficient, and $T_k$ and $T_{k-1}$ are the temperatures at generations $k$ and $k-1$, respectively. We use $P_T = \min\{1, \exp(-\Delta E/T)\}$ as the acceptance probability where $\Delta E$ is the change in objective value, and $T$ is the temperature. To compromise between solution quality and search efficiency, the size of Metropolis sampling $\theta$ is set as the same as the number of blocks. Neighborhood functions are defined in the previous section.

We also considered a SA variant, Fixed Tempreture Algorithm (FTA), where the temperature value remains fixed throughout the entire search. FTA can be obtained by omitting the annealing schedule in SA.

### 5.6 Scatter Search (SS)

Scatter Search was introduced to solve integer programming (Glover, 1998) and was adopted as an effective algorithm for scheduling (Nowicki & Smutnicki, 2006), including for distributed scheduling (Naderi & Ruiz, 2014; Yang, Li, Wang, Liu, & Luo, 2017). Glover (1998) identified a template for Scatter Search, as shown in Fig. 8.

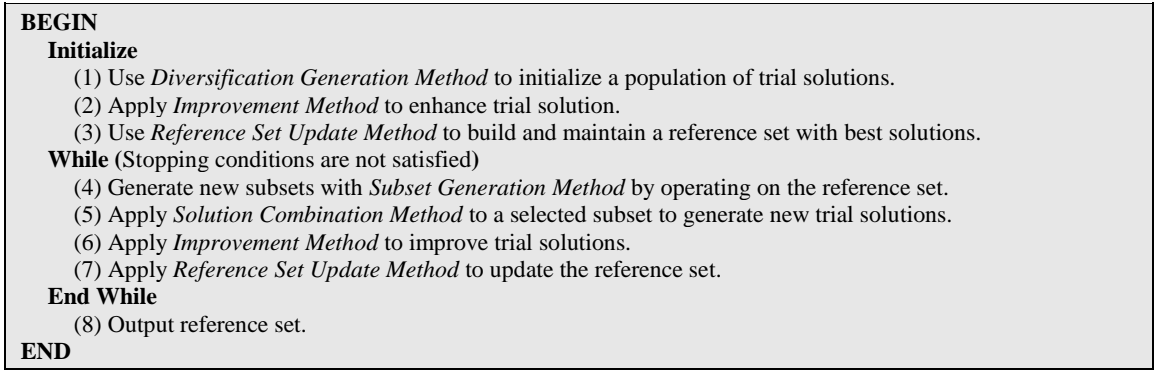

**Fig. 8** Framework for Scatter Search

We sophisticatedly designed all components prerequisite for instantiation of a Scatter Search. First, to initialize trial solutions with high quality and rich diversity, a new block-based NEH (NEH-B) strategy was proposed for *Diversification Generation*. Second, to guarantee diversity in breeding, *Subset Generation Method* adopted the Job Precedence Rule-based distance (Liu, Xu, Qian, Wang, & Chu, 2013) to find the solution furthest away from the best one. Third, the min-max construction based on votes was adopted in *Solution Combination* to generate a new trial solution that was close to parents and inherited more from the parent with higher fitness value. Fourth, to enrich search behaviors and prevent premature convergence, Simulated Annealing was used as *Improvement Method*. Next, we detail implementations.

5.6.1 Diversification Generation Method: Block-based NEH Heuristic (NEH-B)

Diversification generation method generates an initial population of trial solutions. We develop a block-based Nawaz-Enscore-Ham (NEH) heuristic, as detailed in Fig. 9. By running it multiple times, we obtain multiple diversified high-quality solutions that will be fed to Scatter Search as initial solutions.

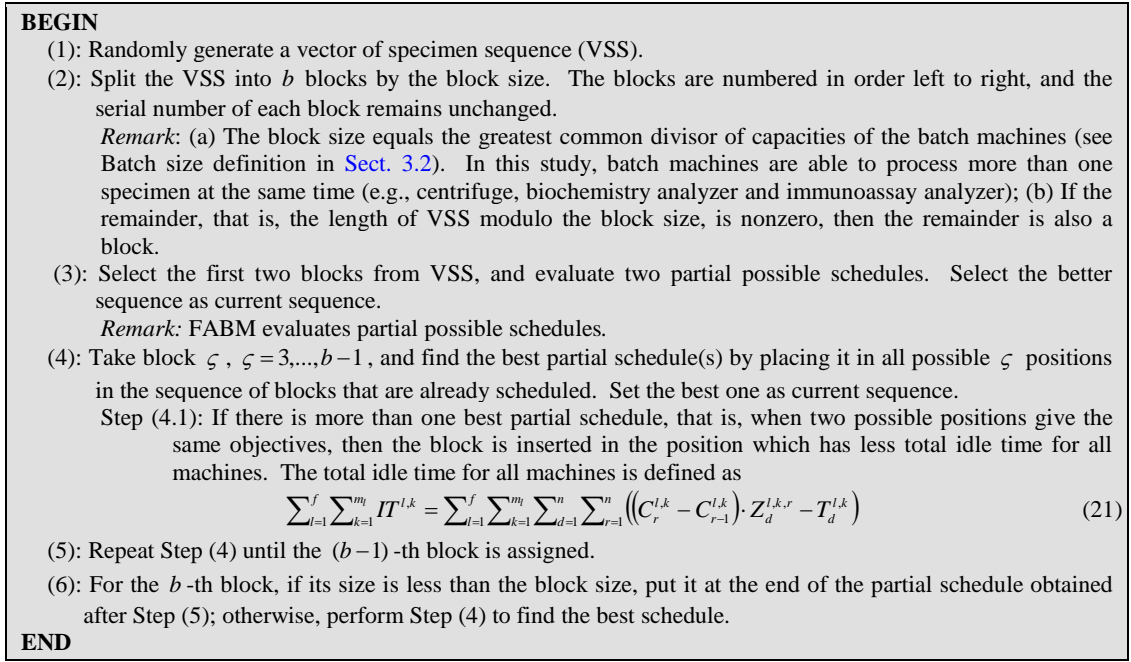

**Fig. 9** Block-based NEH heuristic (NEH-B) for Diversification Generation

NEH heuristic uses a priority rule to gradually build a complete solution, by positing that jobs with higher total processing time should be given higher priority (Nawaz, Enscore, & Ham, 1983). Due to its simplicity and efficacy, NEH has been listed among the best heuristics for flow shops, and is used to construct high-quality solutions fed as initial solutions to metaheuristics (Liu, Wang, Liu, & Wang, 2011).

To make NEH heuristic suitable for our case, and to balance diversification and intensification in search, two enhancements are made to construct the NEH-B heuristic. First, considering alternation of specimen's priority within a batch cannot improve the optimization criterion. We split testing sequence into blocks, regard each block as an entity, and perform NEH operation on these entities rather than on single specimens, avoiding useless moves within the entity. Computational cost is reduced from $O(mn^2)$ in NEH to $O(mb^2)$ in NEH-B, where $m$, $n$ and $b$ denote number of machines, number of specimens and number of blocks, respectively, with $b<n$. Second, we adopt the idea of Ribas et al. (2010) which is: If two possible positions give the same criterion, job is inserted in the position that has less total idle time for all machines. It could save extra cost induced from idle machines and potentially produce solutions with high quality.

### 5.6.2 Improvement: Simulated Annealing

Improvement method performs fine search by exploiting the promising region around the trial solution with a goal of achieving enhanced search quality. We use the simulated annealing as the Improvement method.

### 5.6.3 Subset Generation

Subset Generation generates a new subset by operating on reference set. Our method generates a two-element subset comprising the best solution and the solution furthest away from the best one. The inter-solution distance is measured by Job Precedence Rule-based distance, JPR distance (Liu et al., 2013). JPR distance metric is detailed in Appendix D. The guaranteed diversity benefits breeding diversified solutions for the next Solution Combination procedure.

### 5.6.4 Solution Combination: Min-Max Construction based on Votes

Solution Combination generates a new trial solution by splicing the reference solutions in subset. We adopt min-max construction based on votes by Glover (1994), as shown in Fig. 10.

**BEGIN**

(1): Initialize vote vector $v^t = [v_1^t, v_2^t]'$ as $[0, 0]'$ at $t = 0$, where $v_1^t$ and $v_2^t$ are votes for reference solutions $\pi_1$ and $\pi_2$ for the $t$-th position in trial solution (from left to right, $t = 1, 2, ..., n$). Assign $[MTAT(\pi_2), MTAT(\pi_1)]' ./ (MTAT(\pi_1) + MTAT(\pi_2))$ to influence weight $w = [w_1, w_2]'$ where $MTAT(\pi_1)$ and $MTAT(\pi_2)$ are objective values for $\pi_1$ and $\pi_2$, respectively.

(2): For each reference solution, from left to right, select the first element not included in trial solution as a candidate.

(3): Fill the first unassigned position left to right in trial solution by the min-max construction based on votes.

   (3.1): If two candidates are identical, put any element into the $t$-th position of trial solution, update vote as $v^t = v^{t-1}$, and go to Step (4).

   (3.2): If two candidates are different, select one with respect to minimizing

$$\arg\min_{j \in \{1, 2\}} \| (v^{t-1} + I_j) ./ \sum_{i=1}^{2}(v_i^{t-1} + I_j\{i\}) - w \|_1 \quad (22)$$

where $I_1 = [1, 0]'$, $I_2 = [0, 1]'$, and $\|\cdot\|_1$ is $l^1$-norm, meaning the summation of abstract values.

      (3.2.1): If both $I_j$, $j \in \{1, 2\}$, lead to the same departure in Eq. (22), then randomly select a reference solution's candidate. Update vote as $v^t = v^{t-1} + I_j$ when $j$-th reference solution is selected, and go to Step (4).

      (3.2.2): Select the candidate from the $j$-th reference solution, leading to the minimum departure in Eq. (22), and update vote as $v^t = v^{t-1} + I_j$.

(4): Repeat Steps (2) – (3) until all positions in trial solution are filled.

**END**

**Fig. 10** Min-max construction based on votes for Solution Combination

Its basic idea is to minimize the maximum deviation from the trial solution to the two reference solutions. Two vectors are defined, influence weight vector and vote vector. Influence weight, reflecting reference solution's influence on the trial solution, takes reference solution's fitness value. Here, fitness value is the reciprocal of the objective value. The higher the fitness value of the reference solution, the greater its influence on the trial solution. Vote is defined as accumulated contribution of elements from a reference solution to the trial solution. The more elements the reference solution contributes to the trial solution, the greater the vote. Min-max construction based on votes optimally selects element from reference solutions to construct the trial solution with respect to minimizing the departure of votes (actual contribution) from weights (expected influence): see optimization problem (22) in Fig. 10. $l^1$-norm is used as the distance for measuring the departure of votes from weights. The trial solution is constructed in a position-by-position way in which the first unassigned position from left to right is filled by an element chosen from one of the reference solutions. Votes are updated in the way that the reference solution whose element has been selected would be accredited with one vote. Repeat the procedure until all positions in trial solution are filled with elements.

We illustrate Solution Combination using a concise example in Fig. 11. Two reference solutions $\pi_1$ and $\pi_2$ are $[3, 1, 6, 4, 5, 2]$ and $[3, 4, 6, 1, 5, 2]$, and their objective values $MTAT(\pi_1)$ and $MTAT(\pi_2)$ are 1569.50 and 1829.17, respectively according to FABM. Influence weight is $w = [w_1, w_2]' = [1829.17, 1569.50]' ./ (1569.50 + 1829.17) = [0.54, 0.46]'$, and vote is $v^0 = [0, 0]'$.

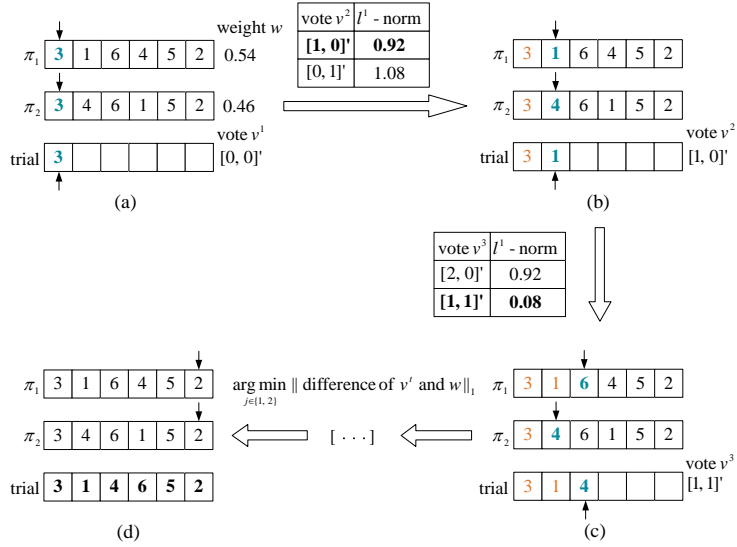

**Fig. 11** Illustration of Solution Combination Method

For each reference solution, at the beginning, the first elements that could fill the first position in trial solution are 3. By Step (3.1), element 3 enters into the first position in trial solution, and vote vector is updated as $v^1 = [0, 0]'$, as shown in Fig. 11a. Then, element 1 in $\pi_1$ and 4 in $\pi_2$ are the first elements not included in trial solution, and are candidates for the second position in trial solution. According to Step (3.2), 1 is selected from $\pi_1$, as shown in Fig. 11b. To be specific, when 1 in $\pi_1$ is selected, vote vector is updated as $v^2 = v^1 + I_1 = [1, 0]'$, and the difference between votes and weights in $l^1$-norm is 0.92; and if 4 in $\pi_2$ was selected, the vote vector was updated as $v^2 = v^1 + I_2 = [0, 1]'$, and the difference was 1.08. Thus, 1 is selected from $\pi_1$ to enter into the second position of trial solution. For the third position in trial solution, 6 in $\pi_1$ and 4 in $\pi_2$ are the first elements not ever having been included in trial solution. A choice of 4 in $\pi_2$ leads to a minimum departure from votes to weights, as shown in Fig. 11c. To be specific, when 6 in $\pi_1$ was selected, vote vector was updated as $v^3 = v^2 + I_1 = [2, 0]'$, and the difference between votes and weights in $l^1$-norm was 0.92; and if 4 in $\pi_2$ is selected, vote vector is updated as $v^3 = v^2 + I_2 = [1, 1]'$, and the difference is 0.08. Repeat the selection until all positions in trial solution are filled with elements. The final combined solution is shown in Fig. 11d.

5.6.5 Reference Set Update

Reference set update method builds and maintains a reference set consisting of best solutions found so far. The improved trial solution from Improvement enters into the reference set, and the solution furthest away from the best one in the reference set is removed. Update the best solution once the improved trial solution outperforms the current best one.

## 6. Tune-ups of the Syntactic Adaptive Problem Solver by Landscape Analysis

In this section, we learn the syntactic adaptive problem solver following the training procedure proposed in Sect. 5. Besides designing a large-scale benchmark, we measure inter-neighbor distances, reveal how landscapes change under different neighborhoods through three landscape analyzers, and determine landscape characteristics that highly influence search performance. We speculate about a chain relatedness starting from search strategy's competence (measured by inter-neighbor distance), then landscape structure, and ending at search performance. Consequently, neighborhood function choice could be inferred as the one that makes search easier.

**6.1 Benchmark**

We design a benchmark that contains 20 toy instances (2, 4, 6, and 8 specimens) and 150 realistic instances (100, 200, 300, 400, and 500 specimens). In our field investigation in a clinical laboratory attached to China-Japan Friendship Hospital in Beijing, average number of tubes per shift is around 200, but for deeper analysis, we enlarged tube number to 500. Thus, our benchmark reflects normal and extreme scenarios in scheduling of medical tests in a clinical laboratory, resulting in a more discriminating benchmark.

To generate toy instances, several modifications are made in settings in Table 1. The number of centrifuge machines in Region 2 is changed to 2, and capacities for centrifuge, biochemical analyzer, and immunoassay analyzer are all set to be 2. The ratio of biochemical tests to immunologic tests is one to one. For each test size (2, 4, 6, 8), five instances are generated, generating a total of 20 toy instances. To generate realistic instances, the settings in Table 1 are adopted. Ratios of biochemical tests to immunologic tests are 1:4, 1:1, and 4:1, respectively. For each combination of test size (100, 200, 300, 400, 500) and each test ratio, ten instances are generated, for a total of 150 instances. The term "*INSTANCE_BT_IT_Idx*" denotes *BT* specimens for biochemical tests, *IT* specimens for immunologic tests, and *Idx* as the index number for each combination of *BT* and *IT*. Benchmark is summarized in Table 3. Instances are publicly available at Wang and Liu (2022). All experiments were performed on an Intel Xeon Gold 6254 3.10GHz machine with 128 GB memory running Ubuntu 18.04.3 LTS. We experimented with Gurobi v9.5.0 and MATLAB 2020b.

**Table 3** Summary of the benchmark

| Instances | Number of specimens | Ratio of biochemical tests to immunologic tests | Number of instances per combination | Total instances |
|---|---|---|---|---|
| Toy instance | {2, 4, 6, 8} | {1:1} | 5 | 20 |
| Medium-size | {100} | {1:4, 1:1, 4:1} | 10 | 30 |
| Large-size | {200, 300, 400, 500} | {1:4, 1:1, 4:1} | 10 | 120 |

**6.2 Inter-neighbor Distance under Different Neighborhoods**

Inter-neighbor distance distinguishes the ability of neighborhood's exploitation over search space. We theoretically measure inter-neighbor distance under each neighborhood. To save space, proofs are put in Appendix A. Results are given in Table 4.

We conduct experiments on representative instances and verify that experimental values are consistent with theoretical values; see Table 5. For medium- and large-size instances, Inverse neighborhood produces the largest inter-neighbor distance, suggesting a potential largest exploration over the solution space. Swap and Insert have the smallest inter-neighbor distances, favoring smaller exploitation. Block-based Insert is at a moderate perturbation level among the four neighborhoods.

**Table 4** Theoretical expectations and variances of the inter-neighbor distance under different neighborhoods

| Neighborhood | Expectation of Inter-neighbor Distance | Variance of Inter-neighbor Distance |
|---|---|---|
| Inverse | $(n-2)(n-3)/(6n(n-1))$ | $(n+1)(n-2)(n-3)(7n-18)/(180n^2(n-1)^2)$ |
| Insert | $2(n+1)/(3n(n-1))$ | $2(n+1)(n-2)/(9n^2(n-1)^2)$ |

| Swap | $2(2n-1)/(3n(n-1))$ | $8(n+1)(n-2)/(9n^2(n-1)^2)$ |
| Block-based Insert | $2n(b+1)/(3b^2(n-1))$ | $2(b+1)(b-2)n^2/(9b^4(n-1)^2)$ |

**Table 5** Theoretical and experimental inter-neighbor distances under different neighborhoods on representative instances

| Instance | Neighborhood | Theoretical values | | Experimental values | |
|---|---|---|---|---|---|
| | | Mean | Variance | Mean | Variance |
| Medium-size instance (INSTANCE_20_80_5) | Inverse | 0.1600 | 0.0371 | 0.1678 | 0.0381 |
| | Insert | 0.0068 | 2.24E-05 | 0.0075 | 2.54E-05 |
| | Swap | 0.0134 | 8.98E-05 | 0.0148 | 1.05E-04 |
| | Block-based Insert | 0.0280 | 3.47E-04 | 0.0259 | 3.24E-04 |
| Large-size instance (INSTANCE_40_160_5) | Inverse | 0.1633 | 0.0380 | 0.1714 | 0.0407 |
| | Insert | 0.0034 | 5.58E-06 | 0.0035 | 5.74E-06 |
| | Swap | 0.0067 | 2.23E-05 | 0.0072 | 2.35E-05 |
| | Block-based Insert | 0.0137 | 8.79E-05 | 0.0132 | 8.60E-05 |
| Large-size instance (INSTANCE_60_240_5) | Inverse | 0.1644 | 0.0383 | 0.1645 | 0.0373 |
| | Insert | 0.0022 | 2.48E-06 | 0.0025 | 2.66E-06 |
| | Swap | 0.0045 | 9.91E-06 | 0.0045 | 9.70E-06 |
| | Block-based Insert | 0.0090 | 3.92E-05 | 0.0087 | 3.82E-05 |
| Large-size instance (INSTANCE_80_320_5) | Inverse | 0.1650 | 0.0384 | 0.1635 | 0.0390 |
| | Insert | 0.0017 | 1.39E-06 | 0.0019 | 1.57E-06 |
| | Swap | 0.0033 | 5.57E-06 | 0.0035 | 5.99E-06 |
| | Block-based Insert | 0.0068 | 2.21E-05 | 0.0069 | 2.29E-05 |
| Large-size instance (INSTANCE_100_400_5) | Inverse | 0.1653 | 0.0385 | 0.1436 | 0.0302 |
| | Insert | 0.0013 | 8.91E-07 | 0.0013 | 7.90E-07 |
| | Swap | 0.0027 | 3.56E-06 | 0.0026 | 3.27E-06 |
| | Block-based Insert | 0.0054 | 1.42E-05 | 0.0054 | 1.40E-05 |

### 6.3 Fitness Distance Correlation Measuring Local-Global Optima Connectivity

We investigate the local-global optima connectivity under each neighborhood function via fitness distance correlation. Two representative instances, INSTANCE_20_80_5 (medium-size) and INSTANCE_60_240_5 (large-size) were selected for illustrative analysis. In the preliminary experiments, Scatter Search was chosen with an iteration number of 20,000. Four algorithms can be derived by specifying one neighborhood in Scatter Search, denoted as SS-SWP (-INS, -INV, or -INB) where Swap (Insert, Inverse, or Block-based Insert) is selected respectively. Objective value relative to the best-so-far solution (*y*-axis) is depicted against its distance from the best-so-far solution (*x*-axis) in Fig. 12.

It is observed that distances between local and global optima (the best-so-far solution) are often less than 0.6 on medium-size instance and 0.5 on large-size instance, indicating that local optima are clustered in a small area of solution space, and local optima of large-size instance are more concentrated. Second, fitness distance correlation under each strategy is between [0.13, 0.47] on the medium-size instance, while in between [0.39, 0.79] on the large-size instance, indicating a stronger connectivity on large-size instance. In particular, on medium-size instance, Inverse and Swap neighborhoods have weaker connectivity than Insert and Block-based Insert neighborhoods. On the large-size instance, except for the Inverse, all neighborhoods have a fitness distance correlation larger than 0.6, featuring a strong connectivity among optima.

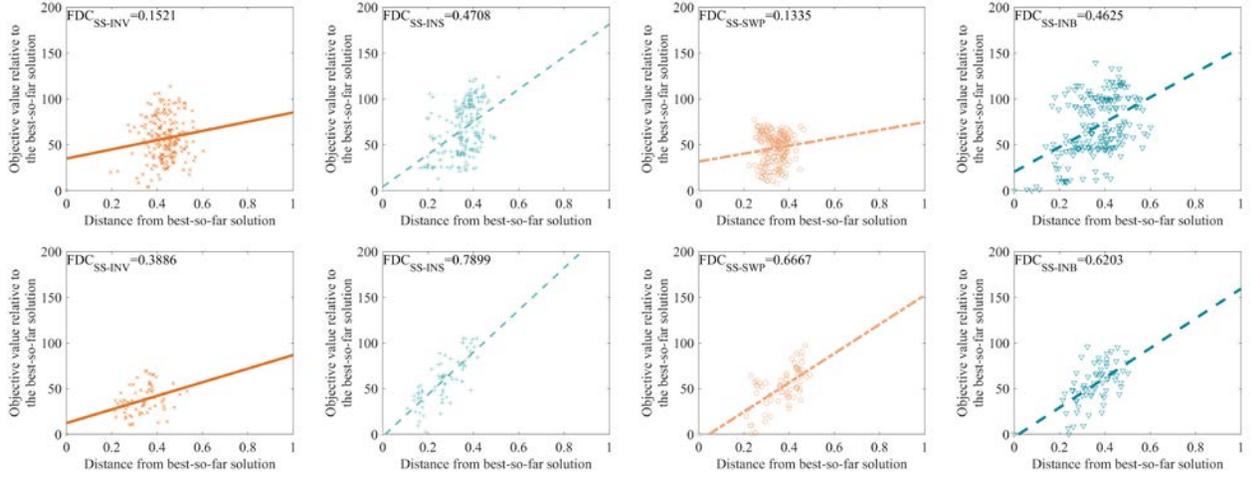

**Fig. 12** Fitness distance correlation measuring local-global optima connectivity. Top: medium-size instance (INSTANCE_20_80_5); Bottom: large-size instance (INSTANCE_60_240_5)

### 6.4 Autocorrelation Coefficient Measuring Landscape's Ruggedness

We utilize autocorrelation coefficient to reveal landscape's ruggedness. Two representative instances are solved by Scatter Search with number of iterations set as 20,000. Autocorrelation coefficient under each neighborhood function is calculated according to Eq. (20). Descending curves and autocorrelation coefficient are depicted for medium-size instance in Fig. 13a and for large-size instance in Fig. 13b.

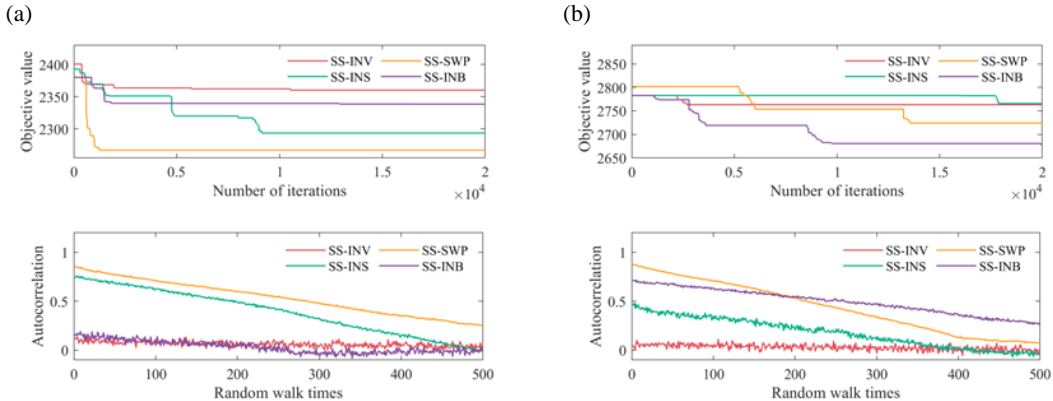

**Fig. 13** Autocorrelation coefficient to measure ruggedness of landscape. (a) medium-size instance (INSTANCE_20_80_5). Top: Descending curves; Bottom: Autocorrelation coefficient along a random walk of 500 steps. (b) large-size instance (INSTANCE_60_240_5)

In Fig. 13a, SS-SWP (-INS) have larger autocorrelation coefficients than SS-INV (-INB), suggesting Swap and Insert neighborhoods lead to smoother landscapes which support a continuously improved objective value, whereas Inverse and Block-based Insert lead to more rugged landscapes. Strategies perform well on a smooth landscape and a fast convergence is expected. This is validated by faster descending curves under SS-SWP (-INS). Besides, as revealed in Table 5, Swap and Insert generate smaller perturbations to solutions than Inverse and Block-based Insert. It can be concluded that medium-size instance is easier under Swap and Insert than Inverse and Block-based Insert since Swap and Insert generated smaller perturbations and led to a relatively smooth landscape.

In Fig. 13b, SS-INB (-SWP) have larger autocorrelation coefficients than SS-INV (-INS), indicating Swap and Block-based Insert neighborhoods lead to smoother landscape. SS-INB achieved the fastest convergence, followed by SS-SWP). As revealed in Table 5, compared with other strategies Block-based Insert generates moderate perturbations to solutions in terms of mean inter-neighbor distance from which we concluded that large-size instance is easier for Block-based Insert that leads to smooth landscapes.

**6.5 Local Optima Network Revealing Difficulty in Jumping out of Local Optima**

Local optima networks under SS-INV (-INS, -SWP, -INB) are visualized in Fig. 14 for toy instance, and in Fig. 15 for medium-size instance. The heat colors palette, a sequential color scheme skewed to reds and yellows, is used. Red identifies the funnel with the best-so-far optima, and a yellow color gradient reflects decrease in solution performance.

a. SS-INV

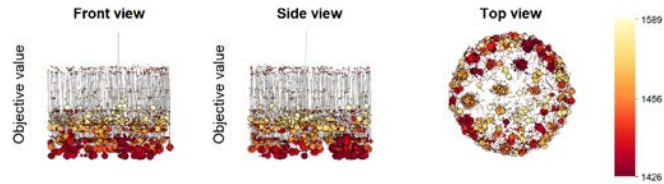

b. SS-INS

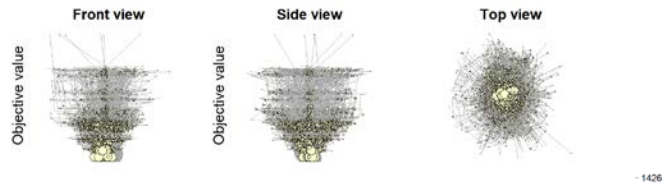

c. SS-SWP

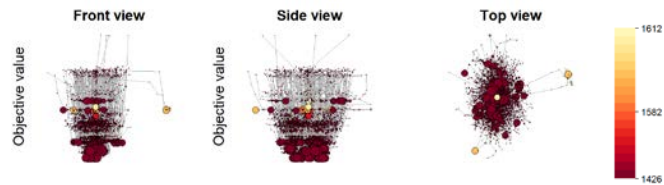

d. SS-INB

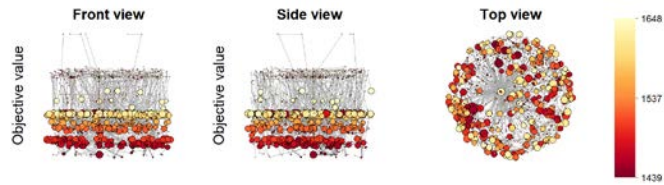

**Fig. 14** Local optima networks under SS on toy instance (INSTANCE_4_4_5)

Since plateaus harass search, it is visually observed from Fig. 14 that even toy instance's local optima networks contain multiple sub-optimal plateaus (attractor plateaus), implying the search complexity for distributed scheduling problems. In Sect. 7.3 we show toy instances' hardness by comparative investigations between SS and Gurobi Optimizer.

For the toy instance in Fig. 14, SS-INS (-SWP) discover escape edges from plateaus to promising regions more easily than SS-INV (-INB) did. It is observed that local optima networks under SS-INV (-INB) are structured as multiple funnels. The considerable amount of attractor plateaus are large-sized, exhibiting great gravitational force, which make searches less likely to escape them. Additionally, SS-INV (-INB) have more scattered sink plateaus, meaning their searching behaviors at the end of search remain diversified, impeding convergence. It is validated by a relatively poor solution under SS-INV (-INB) in Sect. 7.3. It allows the conclusion that toy instances are easier under Swap and Insert (than under Block-based Insert and Inverse) that lead to local optima networks with less funnels and smaller-sized attractor plateaus.

For a medium-size instance in Fig. 15, compared with SS-INV (-INB), SS-INS (-SWP) have smaller plateaus; see plateau average size in Table 6. Smaller average size of plateaus in a local optima network suggests smaller basins of attraction which impose fewer impediments on jumping out of plateaus and support a better exploration over solution space. SS-INS (-SWP) could discover escape edges from attractor plateaus to more promising regions. Search traces under SS-SWP are restricted in small search areas as shown in Fig. 15c, while search traces under SS-INV (-INS, -INB) are more disperse. In this regard, SS-SWP showed aggressive exploitation leading to the fastest convergence.

a. SS-INV

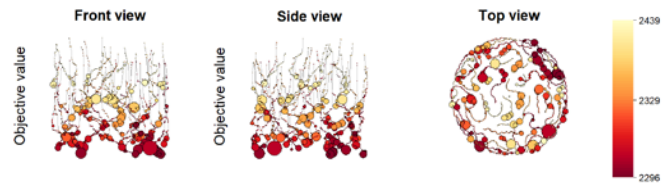

b. SS-INS

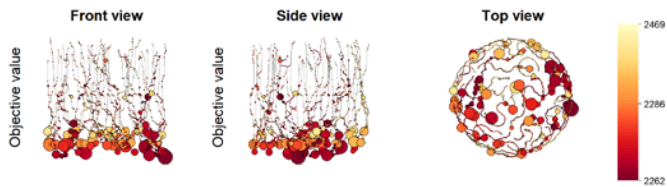

c. SS-SWP

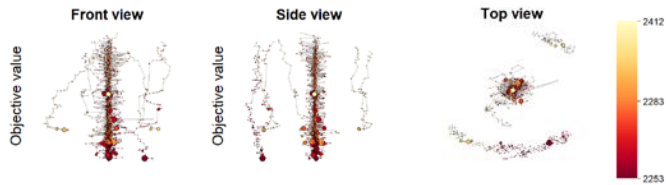

d. SS-INB

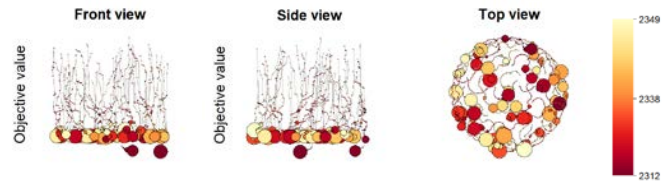

**Fig. 15** Local optima networks under SS on medium-size instance (INSTANCE_20_80_5)

For a large-size instance, since under different neighborhood functions it becomes hard to visually discern differences among local optima networks, statistics abstracted from these networks are used to speculate landscape structures. Compared with SS-INV (-INS), SS-INB (-SWP) have smaller plateaus (see plateau average size in Table 6), indicating less attraction from attractor plateaus and facilitating escape from them. This agrees with observations in autocorrelation coefficients in Sect. 6.4, where SS-INV (-INS) have smaller autocorrelation coefficients, leading to more rugged landscapes and slower convergence rate. Be-cause, as shown in Table 5, compared with Insert and Inverse, Block-based Insert generates moderate perturbations in search space, it can be concluded that large-size instances are easier under moderate-disturbance neighborhood, e.g., Block-based Insert.

**6.6 Strategy Inferred from Landscape Analysis**

6.6.1 Several Interesting Findings Revealed from Preliminary Experiments

In the previous sections, the relationship between domain's syntactic characteristics and neighborhood functions was learned based on preliminary experiments using Scatter Search on representative instances. These experiments revealed several interesting findings, as summarized in Table 6.

**Table 6** Relationship between inter-neighbor distances, landscapes, and performances

| Instance | Neighborhood | Mean Inter-neighbor Distance | Landscape Analyzers | | | Performance Rank |
| --- | --- | --- | --- | --- | --- | --- |
| | | | Fitness Distance Correlation | Autocorrelation AC(1) | Plateau average size | |
| INSTANCE_20_80_5 (100 tasks) | Insert | 0.0068 | 0.4708 | 0.7473 | 27.60 | 2 |
| | Swap | 0.0134 | 0.1335 | 0.8548 | 17.34 | 1 |
| | Inverse | 0.1600 | 0.1521 | 0.1255 | 36.42 | 4 |
| | Block-based Insert | 0.0280 | 0.4625 | 0.1462 | 46.12 | 3 |
| INSTANCE_40_160_5 (200 tasks) | Insert | 0.0034 | 0.8286 | 0.5142 | 23.34 | 4 |
| | Swap | 0.0067 | 0.7985 | 0.5973 | 16.81 | 2 |
| | Inverse | 0.1633 | 0.1813 | 0.1112 | 26.91 | 3 |
| | Block-based Insert | 0.0137 | 0.7639 | 0.5411 | 23.04 | 1 |
| INSTANCE_60_240_5 (300 tasks) | Insert | 0.0022 | 0.7899 | 0.4602 | 22.67 | 4 |
| | Swap | 0.0045 | 0.6667 | 0.8753 | 18.40 | 2 |
| | Inverse | 0.1644 | 0.3886 | 0.0487 | 19.22 | 3 |
| | Block-based Insert | 0.0090 | 0.6203 | 0.7101 | 16.91 | 1 |
| INSTANCE_80_320_5 (400 tasks) | Insert | 0.0017 | 0.7222 | 0.2809 | 13.79 | 4 |
| | Swap | 0.0033 | 0.8023 | 0.8886 | 10.89 | 2 |
| | Inverse | 0.1650 | 0.5334 | 0.0157 | 13.53 | 3 |
| | Block-based Insert | 0.0068 | 0.5190 | 0.7215 | 17.32 | 1 |
| INSTANCE_100_400_5 (500 tasks) | Insert | 0.0013 | 0.4095 | 0.9171 | 18.18 | 4 |
| | Swap | 0.0027 | 0.8170 | 0.9775 | 12.65 | 2 |
| | Inverse | 0.1653 | 0.3640 | 0.0234 | 15.35 | 3 |
| | Block-based Insert | 0.0054 | 0.5009 | 0.8599 | 13.74 | 1 |

For medium-size instances, under the algorithmic framework of Scatter Search, Swap neighborhood showed the best performance. It has the lowest fitness distance correlation, suggesting that fine-grained search is more conducive to enhancing search. Swap neighborhood produces smaller disturbances, favoring smaller exploitation over the solution space. Swap neighborhood has the highest autocorrelation coefficient (smoothest landscapes) and the smallest-sized attractor plateaus (fewest impediments on jumping out of plateaus), supporting a fast descending speed. In contrast, Inverse neighborhood performed worst. It also has very low fitness distance correlation, but lacks the ability to perform fine searches due to its largest disturbances to solutions. Inverse neighborhood was faced with the most rugged landscapes (smallest autocorrelation coefficient), including larger attractors (larger plateau) that made the search stagnate.

For large-size instances, Block-based Insert neighborhood achieved the best performance. It has a moderate fitness distance correlation, sustaining that the search in between fine-grained and coarse-grained searches would be preferable. Block-based Insert could generate moderate disturbances, meeting the re-quirement. Under the moderate disturbances, it generates smoother landscapes with smaller plateaus, thereby fostering search. In contrast, Insert neighborhood performed worst. It has a very high fitness distance correlation, supporting that large-disturbance exploration is an obvious need. Nevertheless, Insert neighborhood has the smallest disturbance, which can not fulfill the need. Insert neighborhood can hardly traverse the search space manifested in rugged landscapes with large-sized plateau due to lack of larger momentum.

6.6.2 Verification of the Relationship between Syntactic Features and Neighborhoods on Other Algorithms

In this section, we will investigate whether the learned relationship holds for other combinatorial optimization algorithms.

We conducted large-scale experiments to test the performance of Scatter Search, Simulated Annealing and Fixed Temperature Algorithm using different neighborhoods on all instances. All algorithms used the same termination criterion, that is, the algorithm terminates when the maximum number of function evaluations reached 10,000. Each algorithm was independently repeated 30 times on each instance. We calculated Average Relative Percentage Deviation values (ARPD, see Sect. 7.2 Performance metric) across all instances for different optimization algorithms invoking various neighborhoods, and grouped ARPD values by instance size. As shown in Fig. 16, for example, when the instance size is 100, the red dots represent the ARPD values when the instance size is 100 using the optimization method that calls the INV neighborhood. There are 30 instances of each size, corresponding to 30 red dots. Lines with different colors are fitted lines of linear regression for 30*5 observations with instance size as independent variable and ARPD as dependent variable. The gray shading on the line is the 95% confidence interval for linear regression. In this way, we plot in Fig. 16 the performance of each algorithm using different neighborhoods against instance size.

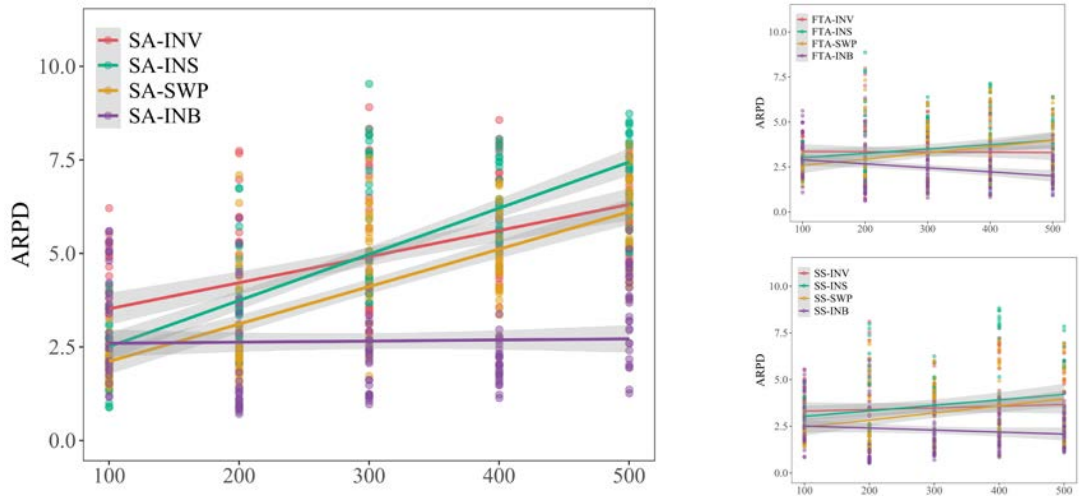

**Fig. 16** Average Relative Percentage Deviation grouped by instance size for each algorithm using different neighborhoods. Different colored lines reflect ARPD of algorithms using different neighborhoods against the instance size. The greater the slope of the line, the greater the performance degradation of algorithms using that neighborhood, and vice versa

It can be observed from Fig. 16 that (*i*) all algorithms using SWP neighborhood achieved the smallest ARPD value on medium-size instances (100 tasks) and performed best. As the instance size increases, the ARPD value keeps increasing, as shown by the yellow lines, indicating that performances of SWP neighborhood gradually deteriorated. (*ii*) On large-size instances (200-500 tasks), all algorithms using the INB neighborhood achieved the smallest ARPD value and performed best. As the instance size increases, all algorithms using the INB neighborhood (purple line) have the smallest slope among all neighborhoods. The purple lines have negative slope values in SS and FTA respectively and are almost flat in SA, indicating that INB neighborhood has the smallest performance degradation among all neighborhoods.

Taken together, through the offline training, we set up an adaptive problem solver based on domain's syntactic characteristics (instance size and local-global optima connectivity, ruggedness, plateaus size) to infer the most suitable neighborhood function. Next, the syntactic adaptive problem solver will be tested, and statistical analysis is performed to validate its performances.

## 7. Experiments and Comparisons

In this section we will test the trained syntactic adaptive problem solver's performances on large-scale instances. It is extensively compared with two constructive heuristics (NEH, NEH-B), nine conventional combinatorial optimization algorithms where neighborhood functions are user-specified, three adaptive solvers based on meta-Lamarckian learning, and Gurobi Optimizer.

**7.1 Algorithms for Comparisons**

7.1.1 Our Proposed Syntactic Adaptive Problem Solvers

As observed from the analysis in the previous section, toy and medium-size instances (100 tasks) were easier under Swap-based neighborhood that generated smaller perturbations and led to smoother landscape with smaller plateaus. For large-size instances (200 to 500 tasks), INB, the new neighborhood that allowed a quick step over multiple funnels and plateaus performed best. Accordingly, the syntactic adaptive problem solver selects the SWP and INB neighborhoods to solve medium- and large-size instances, respectively. For toy and medium-size instances, we have three syntactic adaptive problem solvers, namely APS-SA-SWP (SA adopts SWP-based neighborhood), APS-FTA-SWP and APS-SS-SWP. For large-size instances, we have APS-SA-INB, APS-FTA-INB, and APS-SS-INB.

7.1.2 Conventional Combinatorial Optimization Algorithms

Unlike syntactic adaptive problem solvers, neighborhood functions in conventional combinatorial optimization algorithms are specified by the user, rather than derived from learned domain syntactic features. In this study, three combinatorial optimization algorithms, SA, FTA, and SS, and three neighborhood functions (one of the four neighborhoods are used in the syntactic adaptive problem solver), generate a total of nine conventional combinatorial optimization algorithms for comparison. For example, SA-INV stands for SA with Inverse neighborhood.

7.1.3 Generative Adaptive Problem Solvers Leveraging Meta-Lamarckian Learning

For further comparison, we used generative adaptive problem solvers where meta-Lamarckian learning was relied upon for adaptively selecting the most-rewarded neighborhood from the composite neighborhood (four neighborhoods) at runtime. Unlike syntactic adaptive problem solvers that leverage the explicit relatedness between landscape structures and neighborhood functions, generative adaptive problem solvers implicitly learn landscape structures from past problem-solving experiences, without exploiting prior explicit knowledge of the problem landscape. We have three generative adaptive problem solvers, namely SA-ML, FTA-ML and SS-ML. For example, SA-ML represents SA using meta-Lamarckian learning. The relationship between the three generative adaptive problem solvers is as follows. By fixing the temperature in SA-ML, FTA-ML can be realized. Using SA-ML as the improvement method for SS, SS-ML can be realized. Next, we mainly give the SA-ML.

1. Meta-Lamarckian learning to manage composite neighborhood (SA-ML)

Diverse neighborhoods lead to double-edged sword effects. While perturbation behaviors avoid stagnancy, a dilemma occurs over how to select the most appropriate neighborhood at runtime to eliminate negative effects of improper neighborhood. Meta-Lamarckian learning (Ong & Keane, 2004) was introduced to select the most promising neighborhood from composite neighborhood by learning historical search traces for a single instance at runtime. In simulated annealing, meta-Lamarckian learning adaptively selects neighborhood from the composite neighborhood, as depicted in Fig. 17. Its input is a trial solution. It performs fine search by exploiting the promising region around the trial solution.

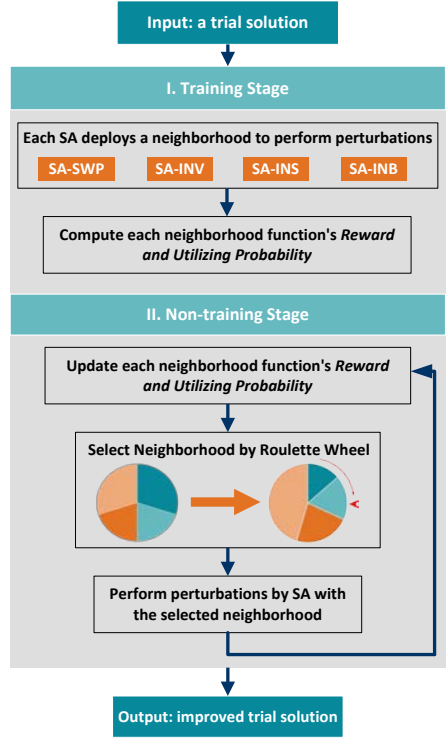

**Fig. 17** Meta-Lamarckian learning enhanced simulated annealing (SA-ML)

Meta-Lamarckian learning, divided into training and non-training stages, has in the training stage four simulated annealing algorithms, each algorithm deploying a neighborhood to perform perturbations with the same size of Metropolis sampling. Reward $\eta_i$ of the $i$-th neighborhood is

$$\eta_i = |MTAT_{b,i} - MTAT_{a,i}|/\theta \tag{23}$$

where $MTAT_{b,i}$ and $MTAT_{a,i}$ are objective values before and after the search using the $i$-th neighborhood, and $\theta$ is the size of Metropolis sampling. To compromise between solution quality and search efficiency, the size of Metropolis sampling is set as the same as the number of blocks.

Utilizing probability for each neighborhood is calculated as

$$p_{ut,i} = \eta_i \Big/ \sum_{j=1}^{K} \eta_j \tag{24}$$

where $p_{ut,i}$ is the utilizing probability for the $i$-th neighborhood, and $K$ is the total number of neighborhoods.

At the non-training stage of meta-Lamarckian learning, based on utilizing probability a roulette wheel rule (Goldberg, 1989) decides which neighborhood to be selected. If the $i$-th neighborhood is selected, its reward will be updated by $\eta_i = \eta_i + \Delta\eta_i$, where $\Delta\eta_i = |MTAT_{b,i} - MTAT_{a,i}|/\theta$ is the reward from the $i$-th neighborhood during the non-training phase. Utilizing probability for each neighborhood is updated again by Eq. (24). Repeat the non-training stage until stopping conditions are satisfied. Interested readers are referred to Ong and Keane (2004).

### 7.2 Performance Metric

For comparisons of algorithms, the following performance metrics are considered.

1. Best relative error to $C^*$ (BRE), where, $C^*$ is the optimal or best-so-far value
2. Average relative error to $C^*$ (ARE)
3. Worst relative error to $C^*$ (WRE)
4. Average Relative Percentage Deviation (ARPD)

The $a$-th algorithm's ARPD on the $k$-th instance is

$$ARPD_{a,k} = \sum_{l=1}^{L}(I(a,k)_l - C_k^*)/C_k^* \cdot L \qquad (25)$$

where $I(a,k)_l$ is the best result on the $l$-th repetition by the $a$-th algorithm on the $k$-th instance, $C_k^*$ is the optimal or best-so-far value for the $k$-th instance, $L$ is the total repetition of each algorithm on each instance.

5. Average CPU time (ACPU)

The $a$-th algorithm's ACPU on the $k$-th instance is

$$ACPU_{a,k} = \sum_{l=1}^{L} CPU(a,k)_l / L \qquad (26)$$

where $CPU(a,k)_l$ is the CPU time in seconds on the $l$-th repetition by the $a$-th algorithm on the $k$-th instance.

### 7.3 Comparison with Exact Method on Toy Instances

We investigated the comparative efficiency between Gurobi Optimizer and the non-exact algorithms in Sect. 7.1 on small-size toy instances. Gurobi MIP solver is called by YALMIP in MATLAB (Lofberg, 2005). The branch & bound algorithm is adopted to solve the MIP. Dual simplex method is used for MIP node relaxations. Presolve parameter is set to the automatic aggressiveness level. TimeLimit for Gurobi runs is set to 3,000 seconds. Terminate the Gurobi solver if it fails to obtain optima within 3,000 seconds. Each non-exact algorithm is terminated if it reaches the maximum iterations 10,000. Each non-exact algorithm is independently run 30 times on each instance. We aggregated the ARPD values (Average Relative Percentage Deviation) and ACPU (average computing time) of all algorithms using the same neighborhood function and showed the performance of different neighborhood functions on the toy instances, as listed in Table 7. Gurobi only obtains optima on INSTANCE_1_1 but fails to find optima on the rest of the instances in 3000 seconds. For instances with larger than 4 specimens, the proposed adaptive problem solvers outperform Gurobi in terms of computational budgets. We noticed in the distributed job shop literature that Naderi and Azab (2014) optimally solved small sized instances up to 12 jobs using CPLEX, another state-of-the-art MIP solver. The distributed scheduling with heterogeneous, flexible job shops in our study is more difficult to solve than the distributed job shop. Next, we will test the performance of the proposed syntactic adaptive solver in complex medium- and large-size instances that are beyond the feasibility of exact methods.

**Table 7** Average Relative Percentage Deviation and Average CPU time for Gurobi Optimizer and non-exact algorithms on toy instances

| Instance | Gurobi | | SWP | | INS | | INB | | INV | | ML | |
|---|---|---|---|---|---|---|---|---|---|---|---|---|
| | **ARPD** | **ACPU** | **ARPD** | **ACPU** | **ARPD** | **ACPU** | **ARPD** | **ACPU** | **ARPD** | **ACPU** | **ARPD** | **ACPU** |
| 1_1_{1-5} | 0 | 1.05 | 0.0629 | 7.07 | 0.0629 | 7.15 | 0.0629 | 7.15 | 0.0629 | 7.28 | 0.0629 | 8.67 |
| 2_2_{1-5} | 0 | 3008 | 0.0505 | 7.56 | 0.0456 | 7.64 | 0.0648 | 7.57 | 0.0751 | 7.64 | 0.0477 | 9.06 |
| 3_3_{1-5} | 0 | 3021 | 0.0461 | 7.78 | 0.0428 | 7.83 | 0.0593 | 7.91 | 0.0789 | 7.81 | 0.0432 | 8.08 |
| 4_4_{1-5} | 0.0037 | 3005 | 0.0087 | 7.98 | 0.0025 | 8.05 | 0.0316 | 8.00 | 0.0151 | 8.13 | 0.0021 | 8.22 |

**7.4 Comparisons on Medium- and Large-Size Instances**

7.4.1 Large-Scale Experimental Results at a Glance

We examine our proposed adaptive problem solvers' performances on 30 medium-size instances and 120 large-size instances. All algorithms are terminated if they reach the maximum iterations 10,000. Each algorithm is independently run 30 times on each instance. In Sect. 6.6, it is learned that SWP is the most effective neighborhood for medium-size instances, and INB is the most effective neighborhood for large-size instances. Results in Table 8 and on-line materials (K. Wang & Liu, 2022) agree with speculations. Large-scale experimental results preliminarily demonstrated that the most efficient neighborhood could be effectively inferred from syntactic features by leveraging landscape structures.

**Table 8** Summary of results on medium- and large-size instances

| Algorithm | Instances(100) | | Algorithm | Instances(200) | | Instances(300) | | Instances(400) | | Instances(500) | |
|---|---|---|---|---|---|---|---|---|---|---|---|
| | **ARPD** | **ACPU** | | **ARPD** | **ACPU** | **ARPD** | **ACPU** | **ARPD** | **ACPU** | **ARPD** | **ACPU** |
| APS-FTA-SWP | 0.0217 | 12.30 | APS-FTA-INB | 0.0214 | 20.53 | 0.0217 | 21.56 | 0.0220 | 26.29 | 0.0333 | 31.80 |
| APS-SA-SWP | **0.0197** | 12.44 | APS-SA-INB | 0.0205 | 21.13 | 0.0202 | 21.38 | **0.0201** | 27.01 | 0.0317 | 30.92 |
| APS-SS-SWP | 0.0211 | 19.55 | APS-SS-INB | **0.0184** | 56.92 | **0.0177** | 124.44 | 0.0211 | 243.24 | **0.0204** | 527.51 |
| FTA-INS | 0.0216 | 12.34 | FTA-INS | 0.0369 | 20.53 | 0.0558 | 21.50 | 0.0647 | 26.24 | 0.0691 | 31.83 |
| SA-INS | 0.0219 | 12.44 | SA-INS | 0.0353 | 20.40 | 0.0540 | 21.31 | 0.0626 | 25.92 | 0.0671 | **30.22** |
| SS-INS | 0.0277 | 19.67 | SS-INS | 0.0327 | 56.78 | 0.0339 | 124.25 | 0.0458 | 243.18 | 0.0333 | 527.26 |
| FTA-INV | 0.0367 | 12.31 | FTA-INV | 0.0384 | 20.43 | 0.0507 | 21.50 | 0.0594 | 26.14 | 0.0588 | 31.66 |
| SA-INV | 0.0357 | 12.42 | SA-INV | 0.0375 | 20.43 | 0.0495 | 21.28 | 0.0574 | **25.89** | 0.0578 | 30.34 |
| SS-INV | 0.0325 | 19.64 | SS-INV | 0.0320 | 56.84 | 0.0319 | 124.35 | 0.0397 | 242.99 | 0.0304 | 527.99 |
| FTA-INB | 0.0324 | 12.36 | FTA-SWP | 0.0313 | 20.40 | 0.0445 | 21.51 | 0.0511 | 26.16 | 0.0577 | 31.68 |
| SA-INB | 0.0327 | 12.80 | SA-SWP | 0.0291 | 20.33 | 0.0432 | **21.25** | 0.0487 | 25.98 | 0.0571 | 30.27 |
| SS-INB | 0.0297 | 19.54 | SS-SWP | 0.0283 | 56.77 | 0.0324 | 124.43 | 0.0391 | 243.36 | 0.0319 | 527.65 |
| FTA-ML | 0.0241 | 12.48 | FTA-ML | 0.0305 | **16.69** | 0.0376 | 22.06 | 0.0367 | 27.31 | 0.0508 | 32.18 |
| SA-ML | 0.0213 | 12.49 | SA-ML | 0.0294 | 16.71 | 0.0349 | 22.09 | 0.0339 | 27.33 | 0.0470 | 32.04 |
| SS-ML | 0.0256 | 19.62 | SS-ML | 0.0250 | 53.18 | 0.0242 | 124.85 | 0.0288 | 244.33 | 0.0243 | 529.40 |
| NEH | 0.0414 | 5.30 | NEH | 0.0545 | 27.34 | 0.0934 | 76.21 | 0.1426 | 164.12 | 0.0687 | 286.18 |
| NEH-B | 0.0414 | **4.27** | NEH-B | 0.0415 | 21.78 | 0.0384 | 61.59 | 0.0427 | 130.46 | 0.0324 | 246.03 |

7.4.2 Impact of Different Neighborhood Functions on Search Performance

To verify the effectiveness of neighborhood selection by syntactic features, we further compared the impact of different neighborhood functions on search performance.

We aggregated the ARPD of all algorithms using the same neighborhood function and showed the performance of different neighborhood functions on medium-size (Fig. 18a) and large-size instances (Fig. 18b), respectively. The INV on the horizontal axis represents the aggregated ARPD using the three algorithms having Inverse neighborhood (i.e., FTA-INV, SA-INV, and SS-INV). SWP in Fig. 18a and INB in Fig. 18b represented the aggregated ARPD of the algorithm using neighborhood selected by syntactic characteristics, respectively. For easy identification, we labeled the Adaptive Problem Solver below SWP and INB. It can be observed from Fig. 18a and Fig. 18b that the syntactic adaptive problem solver performed best on medium-size instances and large-size instances. Specifically, on medium-size instances, using two-sided Wilcoxon's tests, the algorithm with the SWP neighborhood outperformed the algorithm with the INV ($P = 3.64 \times 10^{-16}$), INS ($P = 4.07 \times 10^{-5}$), INB ($P = 2.64 \times 10^{-11}$), and ML ($P = 2.46 \times 10^{-5}$) neighborhoods significantly. On large-size instances, the algorithm with the INB neighborhood performed significantly better than the algorithm with the INV ($P = 1.03 \times 10^{-60}$), INS ($P = 1.46 \times 10^{-60}$), SWP ($P = 2.06 \times 10^{-60}$), and ML ($P = 1.63 \times 10^{-60}$) neighborhoods.

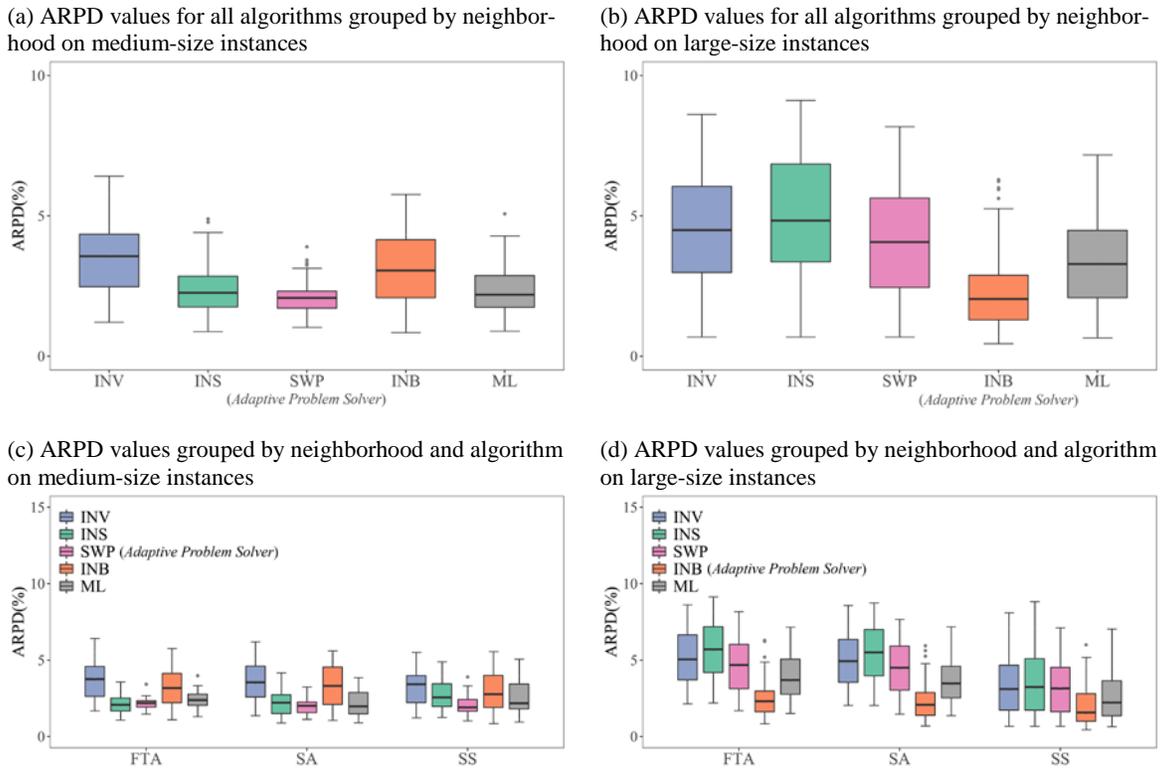

(a) ARPD values for all algorithms grouped by neighborhood on medium-size instances

(b) ARPD values for all algorithms grouped by neighborhood on large-size instances

(c) ARPD values grouped by neighborhood and algorithm on medium-size instances

(d) ARPD values grouped by neighborhood and algorithm on large-size instances

**Fig. 18** Impact of different neighborhood functions on search performances

We performed Friedman's test to identify the existence of differences in average performance among neighborhoods. The null hypothesis $H_0$ is that there is no significant difference among neighborhoods, while the alternative hypothesis $H_1$ indicates the presence of significant differences. Table 9 shows results of Friedman's test, including average rank for each neighborhood and $p$ value. The null hypothesis was rejected under the significance level of 0.001, suggesting that significant differences among neighborhoods exist. Neighborhoods inferred using syntactic features, i.e., SWP and INB, performed best on medium- and large-size instances, respectively, followed by ML.

**Table 9** Results of Friedman's Test

| Instance size | INV | INS | SWP | INB | ML |
|---|---|---|---|---|---|
| 100 | 4.67 | 2.43 | **1.92** | 3.67 | 2.31 |
| 200 | 4.37 | 4.37 | 2.79 | **1.04** | 2.43 |
| 300 | 3.73 | 4.67 | 3.56 | **1.00** | 2.04 |
| 400 | 4.04 | 4.80 | 3.16 | **1.00** | 2.00 |
| 500 | 3.62 | 4.59 | 3.49 | **1.11** | 2.19 |
| 200 to 500 | 3.94 | 4.61 | 3.25 | **1.04** | 2.17 |

The ML, composite neighborhoods managed by meta-Lamarckian learning, achieved competitive performances, implying reliability when prior knowledge of landscape is unknown. Meta-Lamarckian learning's impacts on performances for INSTANCE_60_240_5 are visually investigated in Fig. 19. After the training stage, Block-based Insert is rewarded by the highest utilization probability, yet it is not always chosen along the search. Convergence is faster as Block-based Insert is utilized, while convergence is slower, even kept stagnation as Block-based Insert is not applied.

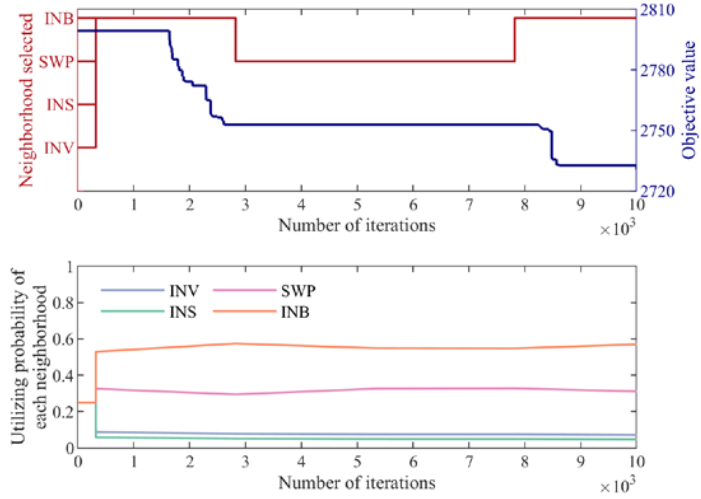

**Fig. 19** Illustrative example on meta-Lamarckian learning's impacts on performances for INSTANCE_60_240_5. Top: Descending curves under SS-ML and actually utilized neighborhoods along search under meta-Lamarckian learning; Bottom: Utilization probability of each neighborhood

We further observe the effect of using different neighborhood functions in an algorithm on search performance. The difference in search ability of each algorithm using different neighborhoods is shown in Fig. 18c (medium-size instances) and Fig. 18d (large-size instances). It can be observed that when the neighborhood specified by syntactic features is used in the algorithm, its performance exceeds that of using other neighborhoods. Specifically, on medium-size instances, FTA-SWP significantly outperformed FTA-INV ($P = 2.88 \times 10^{-6}$), FTA-INB ($P = 2.05 \times 10^{-4}$), and FTA-ML ($P = 0.013$). SA-SWP significantly outperformed SA-INV ($P = 2.35 \times 10^{-6}$) and SA-INB ($P = 4.07 \times 10^{-5}$). SS-SWP significantly outperformed SS-INV ($P = 2.88 \times 10^{-6}$), SS-INS ($P = 1.92 \times 10^{-6}$), SS-INB ($P = 2.88 \times 10^{-6}$), and SS-ML ($P = 6.16 \times 10^{-4}$). FTA-SWP was not significantly different from FTA-INS ($P = 0.910$), and SA-SWP was not significantly different from SA-INS ($P = 0.082$) and SA-ML ($P = 0.178$). On large-size instances, FTA-INB significantly outperformed FTA-INV ($P = 1.97 \times 10^{-21}$), FTA-INS ($P = 1.97 \times 10^{-21}$), FTA-SWP ($P = 2.02 \times 10^{-21}$), and FTA-ML ($P = 1.97 \times 10^{-21}$). SA-INB significantly outperformed SA-INV ($P = 1.97 \times 10^{-21}$), SA-INS ($P = 2.07 \times 10^{-21}$), SA-SWP ($P = 2.80 \times 10^{-21}$), and SA-ML ($P = 2.24 \times 10^{-21}$). SS-INB significantly outperformed SS-INV ($P = 2.60 \times 10^{-21}$), SS-INS

($P = 3.78 \times 10^{-21}$), SS-SWP ($P = 2.95 \times 10^{-21}$), and SS-ML ($P = 4.18 \times 10^{-21}$). These experimental and statistical results justified inferring neighborhoods from syntactic features.

7.4.3 Statistical Comparisons of Algorithms of Algorithms' Performances on Instances of Various Sizes

Fig. 20 shows the differences in ARPD between syntactic adaptive problem solvers and generative adaptive problem solvers, conventional algorithms, constructive heuristics on medium- and large-size instances.

a. Difference in ARPD among algorithms on medium-size instances

b. Difference in ARPD among algorithms on large-size instances with 200 specimens

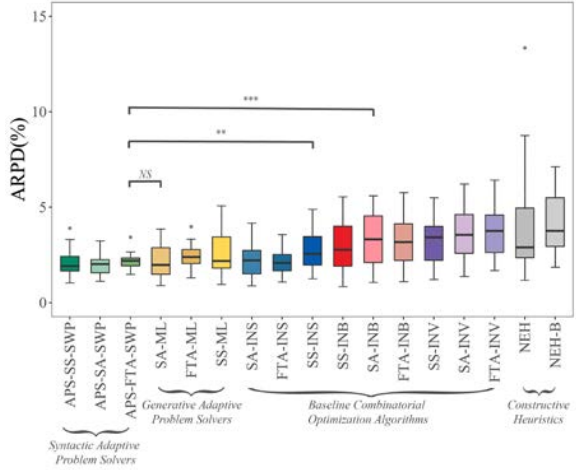
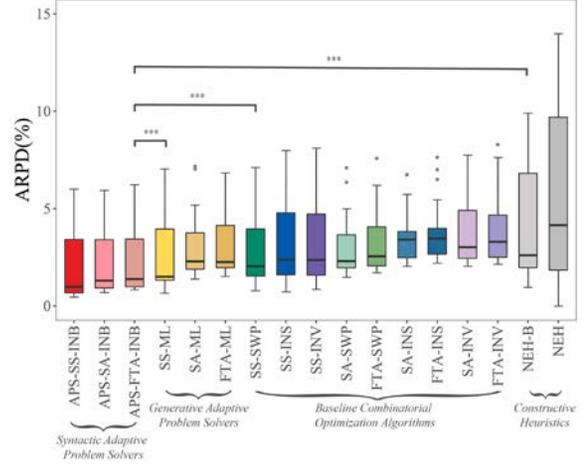

c. Difference in ARPD among algorithms on large-size instances with 300 specimens

d. Difference in ARPD among algorithms on large-size instances with 400 specimens

e. Difference in ARPD among algorithms on large-size instances with 500 specimens

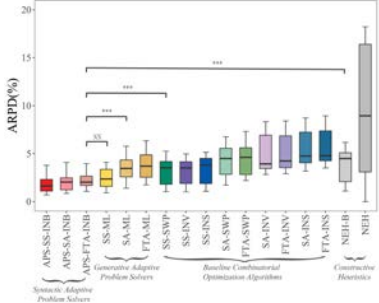
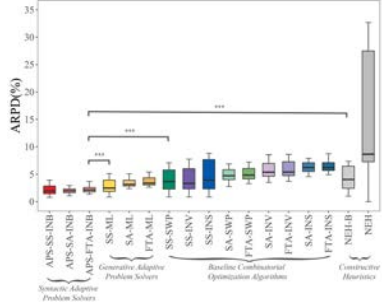
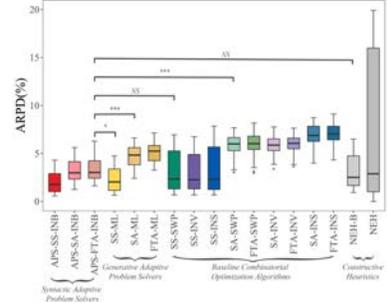

**Fig. 20** Differences in ARPD among syntactic adaptive problem solvers, generative adaptive problem solvers, conventional algorithms, and constructive heuristics on medium- and large-size instances. Statistical comparison is made using two-sided Wilcoxon's tests. *P<0.05; ** P<0.01; ***P<0.001; NS, not significant (P>0.05)

On large instances, taking the 200-specimen instance as an example, as shown in Fig. 20b, we found that (*i*) among the syntactic adaptive solvers, APS-SA-INB ($P = 4.68 \times 10^{-3}$) and APS-SS-INB ($P = 1.36 \times 10^{-5}$) significantly outperformed APS-FTA-INB; (*ii*) Compared with other algorithms, even the worst performing syntactic adaptive solver APS-FTA-INB could significantly outperform all algorithms, e.g., SS-ML ($P = 4.45 \times 10^{-5}$), SA-ML ($P = 1.92 \times 10^{-6}$), FTA-ML ($P = 1.73 \times 10^{-6}$), SS-SWP ($P = 8.47 \times 10^{-6}$), SS-INS ($P = 3.88 \times 10^{-6}$), SS-INV ($P = 3.18 \times 10^{-6}$), SA-SWP ($P = 5.22 \times 10^{-6}$), FTA-SWP ($P = 1.92 \times 10^{-6}$), SA-INS ($P = 2.13 \times 10^{-6}$), FTA-INS ($P = 1.73 \times 10^{-6}$), SA-INV ($P = 1.73 \times 10^{-6}$), FTA-INV ($P = 1.73 \times 10^{-6}$), NEH-B ($P = 1.73 \times 10^{-6}$), and NEH ($P = 8.19 \times 10^{-5}$). We found the same trend on large instances of other sizes, as shown in Figs. 20c through e.

## 8. Conclusions and Future Work

Application of machine learning in clinical medicine still faces daunting challenges. Our research is in the direction of narrowing the gap between machine learning research and clinics-orientated operations practice. This paper formulated a mixed integer programming model for scheduling clinical laboratory testing, learned a syntactic adaptive problem solver by leveraging instances' landscape features, and validated its efficacy in large-scale instances. We made several assumptions, and limitations exist, both of which should be noted.

*Formulation and benchmark instances for real-world clinical laboratory scheduling problem*: We formulated a mixed integer programming model for the heterogeneous, flexible job shop in distributed scheduling (D-HFJSP). This model focused on efficiency criterion in short-term scheduling, specifically deducing turnaround time. To our knowledge, this is the first formulation and benchmark for scheduling related to clinical laboratory and for distributed scheduling in heterogeneous, flexible job shop. Other issues related to efficiency worth considering are tests' priority, instrument recalibration, and staff changeover. In addition to efficiency, clinical laboratory should be robust enough to uncertainties in long-term scheduling, shortage of consumables (reagents and probes), breakdowns and other. Furthermore, a mathematical model is fragile when the situation changes, and incorporation of staff's situational awareness helps to achieve expected-quality schedules. Initiatives related to improved patient care, such as specimen dividing for reduced blood draws, could be considered.

*Adaptive problem solving*: The no-free-lunch (NFL) theorem sustains that algorithm performance is problem-dependent, and no algorithm performs well across all possible problems (Wolpert & Macready, 1997). Awareness of relatedness between landscape structures and search strategies would benefit the identification of an effective strategy. We presented a syntactic adaptive problem solver leveraging instances' landscape characteristics to specify which search strategy to use. The adaptive solver automatically instantiated the best algorithms for various scenarios, reducing repetitive algorithm-tuning efforts as instances frequently change over time. Since our initiative was to provide an easy-to-use solver to the physician, it was essentially a preliminary lookup table. More precisely, the relatedness between characteristics of instance's landscape and algorithm performance was established based on experimental observations but not on statistical analysis. It is worthwhile to further deploy state-of-the-art regression models from machine learning to validate the relationship by enrolling enough observations. Another important improvement for future work is to perform deeper analysis on instance space and to incorporate these features into the syntactic adaptive problem solve. It is worth further bringing in machine learning by acknowledging problems' distribution. It is, however, interesting and challenging to explicitly characterize a distribution.

*Landscape analysis*: Characteristics of landscape are critical indicators of search space topology, determining instance hardness under different neighborhood functions. In our study, it revealed landscape structures, guided strategy selection, and motivated the design of a new strategy. Three landscape analysis analyzers, namely fitness distance correlation, autocorrelation coefficient, and local optimal network were resorted to in the learning domain's landscape characteristics explicitly. We investigated changes in landscape structure and found that disturbance amplitude, local-global optima connectivity, landscape's ruggedness and plateau size fairly predict search strategies' efficacy. The relationship between the features of the landscape under each neighbourhood function and its performance were clearly presented across various instance sizes. These landscape characteristics were recognized as the syntactic characteristics for our proposed adaptive problem solver. It is worth exploiting more landscape analyzers to investigate instance characteristics. It is, however, challenging to achieve a trade-off between solver's performance gains and computational effort required for exploratory landscape analysis.

*Local optima network*: It visualized landscapes in nodes and edges. We used it as an offline approach to successfully reflect instance hardness and predict search performance. To our knowledge, local optima network was used for the first time to analyze a real-world complex scheduling problem. Since each algorithm is independently run multiple times to reveal the complex network, it is too expensive to guide search online. Realizing rapid prototyping of the local optima network could be a possible research direction.

*Meta-Lamarckian learning:* A generative adaptive problem solver was constructed for comparison, where meta-Lamarckian learning was relied upon for implicitly learning the landscape structure. It adaptively selected the most promising strategy by rewards, and achieved beyond average performance. Since at the training stage every strategy needs to run for the same amount of perturbations, consuming a lot of computational budgets, it is interesting to see an economic training strategy to save overhead, especially on expensive-to-evaluation instances. At the training stage, each strategy started the search from the same startpoint at which instance hardness might be different for each strategy, leading to "unfair" rewards. A reward strategy based on both multiple starting points and voting is expected to provide a robust measure.

*Search strategy with adjustable disturbance amplitude*: Changes in landscape were investigated under strategies that imposed different-magnitude perturbations. The relatedness observed motivated the design of a new search strategy with moderate perturbations. Search strategy with adjustable disturbance amplitude is highly needed to perform an expected exploitation over the search space.

Clinical practice will remain dependent on laboratory staff's experience to draw up schedules even though the advance in applied modeling and an adaptive problem solver is vital towards sophisticated, automated scheduling for clinical decision-making. It helps to achieve a better clinical practice, enhance the quality of healthcare, and shape the laboratory of the future.

## Appendix A: Expectation and Variance of Inter-neighbor Distance under Different Neighborhoods

We prove the expectation and variance of inter-neighbor distance under element-based perturbation (Swap, Insert and Inverse) and block-based perturbation (Block-based Insert).

### A.1 Inter-Neighbor Distance under Element-based Perturbation

We assume that any neighbor could be reached by same probability, and define a discrete random variable $\xi$ to represent the number of elements between any two distinct elements in a vector of specimen sequence (VSS). Two distinct elements $J_1$ and $J_2$, randomly selected from the VSS, are located in the $s$-th and $(s+\xi+1)$-th position, respectively. $d(\xi)$ is the inter-neighbor distance, and $p(\xi)$ is the probability distribution function. The mean and variance of $d(\xi)$ can be defined as

$$ED = \sum d(\xi) p(\xi) \tag{A.1}$$

$$VarD = \sum p(\xi)(d(\xi) - ED)^2 . \tag{A.2}$$

A.1.1 Inter-Neighbor Distance under Swap

That $J_1$ and $J_2$ are swapped will change $(2\xi+1)$ element pairs' precedence, leading to

$$d(\xi) = 2(2\xi+1)/(n \cdot (n-1)) \tag{A.3}$$

under the definition of inter-solution distance metric (see Appendix D). For any two randomly selected elements, the probability that the former is $\xi$ positions ahead of the latter one is

$$p(\xi) = 2(n - \xi - 1)/(n \cdot (n-1)) \tag{A.4}$$

where $\xi \in \{0, 1, ..., n-2\}$. We obtain

$$ED_{SWP} = 2(2n-1)/(3n \cdot (n-1)) \tag{A.5}$$

$$VarD_{SWP} = 8(n+1)(n-2)/(9n^2(n-1)^2) \tag{A.6}$$

A.1.2 Inter-Neighbor Distance under Insert

Inserting $J_2$ before $J_1$ leads to

$$d(\xi) = 2(\xi+1)/(n(n-1)) \tag{A.7}$$

and $p(\xi)$ remains the same as before. We obtain

$$ED_{INS} = 2(n+1)/(3n \cdot (n-1)) \tag{A.8}$$

$$VarD_{INS} = 2(n+1)(n-2)/(9n^2(n-1)^2) \tag{A.9}$$

A.1.3 Inter-Neighbor Distance under Inverse

Inverting the subsequence between $J_1$ and $J_2$ leads to

$$d(\xi) = \xi(\xi-1)/(n(n-1)) \tag{A.10}$$

and $p(\xi)$ remains the same as before. We obtain

$$ED_{INV} = (n-2)(n-3)/(6n(n-1)) \tag{A.11}$$

$$VarD_{INV} = (n+1)(n-2)(n-3)(7n-18)/(180n^2(n-1)^2) \tag{A.12}$$

## A.2 Inter-Neighbor Distance under Block-based Insert

The $n$-length sequence is split into $b$ blocks through diversification generation method (Sect. 5.6.1). Each block is composed of $n_c$ elements/specimens. Blocks are numbered in order, left to right, and each block's serial number remains unchanged. Block-based Insert operates on these blocks, randomly choosing from the sequencee two distinct blocks and inserting the latter block before the front block while keeping unchanged specimens' priority within each block. When operating on the sequence, some elements belonging to the same block may be scattered in certain positions in the sequence, resulting in these elements not being connected with other elements in the same block, but the membership of these elements remains constant. We consider two situations, one is that specimens belonging to the same block are not scattered in the sequence, as shown in Fig. 5d of Sect. 5.3; the other is that specimens belonging to the same block are scattered, as shown in Fig. 5e.

### A.2.1 Elements belonging to the Same Block are not Scattered

Two distinct blocks $B_g$ and $B_h$ are randomly selected. $B_g$ is before $B_h$. $\kappa$, a discrete random variable, represents that there are $\kappa$ blocks between $B_g$ and $B_h$. It leads to an inter-neighbor distance

$$d(\kappa) = 2(\kappa+1)n_c^2/(n(n-1)) = 2(\kappa+1)n/(b^2(n-1)), \qquad (A.13)$$

since $(\kappa+1) \cdot n_c^2$ specimen pairs' precedence are changed and a probability distribution function is calculated with

$$p(\kappa) = 2(b-\kappa-1)/(b(b-1)), \qquad (A.14)$$

where $\kappa \in \{0,1,\ldots,b-2\}$.

To obtain

$$ED_{INB} = 2n(b+1)/(3b^2(n-1)), \qquad (A.15)$$
$$VarD_{INB} = 2(b+1)(b-2)n^2/(9b^4(n-1)^2). \qquad (A.16)$$

### A.2.2 Elements belonging to the Same Block are Scattered

Same as the previous definitions, blocks $B_g$ and $B_h$ are two randomly selected blocks, and $B_g$ is before $B_h$. Elements in $B_g$ and $B_h$ are indexed in order as $\{i_{g,1},\ldots,i_{g,q},\ldots,i_{g,n_c}\}$ and $\{i_{h,1},\ldots,i_{h,q},\ldots,i_{h,n_c}\}$, respectively; $i_{g,q}$ is valued as the position of its corresponding element in whole sequence. It leads to an inter-neighbor distance

$$d(g,h) = \sum_{q=1}^{n_c} 2(i_{h,q} - i_{g,1})/(n(n-1)). \qquad (A.17)$$

Probability that both the $g$-th and the $h$-th blocks are selected is

$$p(g,h) = 2/(b(b-1)). \qquad (A.18)$$

Inter-neighbor distance's expectation is

$$ED = \sum_{1 \leq g < h \leq b} \sum_{q=1}^{n_c} \frac{4(i_{h,q} - i_{g,1})}{n(n-1)b(b-1)}. \qquad (A.19)$$

Its second moment is

$$ED^2 = \sum_{h=2}^{b} \sum_{g=1}^{h-1} \frac{1}{(h-1)(b-1)} \left( \sum_{q=1}^{n_c} \frac{i_{h,q} - i_{g,1}}{n(n-1)/2} \right)^2. \qquad (A.20)$$

Its variance is

$$VarD = ED^2 - E^2D. \tag{A.21}$$

Since it is hard to derive explicit expectation and variance, we provide the upper and lower bounds.

*Expectation's upper bound*: Since $i_{g,1} \geq g$ and $\sum_{h=2}^{b}\sum_{q=1}^{n_c} i_{h,q} \leq \sum_{\tau=1}^{n_c(b-1)}(n_c + \tau)$, we obtain

$$\sum_{1 \leq g < h \leq b}\sum_{q=1}^{n_c} i_{h,q} = \sum_{g=1}^{b-1}\sum_{h=g+1}^{b}\sum_{q=1}^{n_c} i_{h,q} \leq \sum_{g=1}^{b-1}\sum_{\tau=1}^{(b-g)n_c}(gn_c + \tau). \tag{A.22}$$

The upper bound is

$$\begin{aligned}
ED &= \sum_{g=1}^{b-1}\sum_{h=g+1}^{b}\sum_{q=1}^{n_c} 4(i_{h,q} - i_{g,1})/(n(n-1)b(b-1)) \\
&\leq \frac{4n_c}{n(n-1)b(b-1)}\sum_{g=1}^{b-1}\left(g(b-g)n_c + \frac{(b-g)((b-g)n_c + 1)}{2} - g(b-g)\right) \\
&= \frac{2n_c}{n(n-1)b(b-1)}\sum_{g=1}^{b-1}\left(b(bn_c + 1) - (2b+1)g + (-n_c + 2)g^2\right) \\
&= \frac{n_c(4bn_c + n_c - 2b + 1)}{3n(n-1)} \\
&= \frac{(4b+1)n + (-2b+1)b}{3b^2(n-1)}
\end{aligned} \tag{A.23}$$

*Expectation's lower bound*: Since $i_{g,1} \leq (g-1)n_c + 1$, we obtain

$$\sum_{1 \leq g < h \leq b}\sum_{q=1}^{n_c} i_{h,q} = \sum_{g=1}^{b-1}\sum_{h=g+1}^{b}\sum_{q=1}^{n_c} i_{h,q} \geq \sum_{g=1}^{b-1}\sum_{\tau=1}^{(b-g)n_c}(g + \tau). \tag{A.24}$$

The lower bound is

$$\begin{aligned}
ED &= \sum_{g=1}^{b-1}\sum_{h=g+1}^{b}\sum_{q=1}^{n_c} 4(i_{h,q} - i_{g,1})/(n(n-1)b(b-1)) \\
&\geq \frac{4n_c}{n(n-1)b(b-1)}\sum_{g=1}^{b-1}\left(g(b-g) + \frac{(b-g)((b-g)n_c + 1)}{2} - ((g-1)n_c + 1)(b-g)\right) \\
&= \frac{2n_c}{n(n-1)b(b-1)}\sum_{g=1}^{b-1}\left((b^2 n_c + 2bn_c - b) + (-4bn_c - 2n_c + 2b + 1)g + (3n_c - 2)g^2\right). \\
&= \frac{n_c(3n_c + 2b - 1)}{3n(n-1)} \\
&= \frac{3n + (2b-1)b}{3b^2(n-1)}
\end{aligned} \tag{A.25}$$

*Variance's upper bound*: The following inequalities hold for any sequence.

$$i_{g,1} \geq g \tag{A.26}$$

$$\sum_{g=1}^{b-1}\sum_{h=g+1}^{b} a_h^2 \leq \sum_{g=1}^{b-1}\sum_{h=g+1}^{b}\left(\sum_{q=1}^{n_c}(n_c(h-1) + q)\right)^2 \tag{A.27}$$

$$\sum_{g=1}^{b-1}\sum_{h=g+1}^{b} a_h \geq \sum_{g=1}^{b-1}\sum_{\tau=1}^{(b-g)n_c}(g + \tau) \tag{A.28}$$

where $a_h = \sum_{q=1}^{n_c} i_{h,q}$.

The upper bound of the second moment is

$$ED^2 = \frac{8}{n^2(n-1)^2 b(b-1)} \sum_{g=1}^{b-1}\sum_{h=g+1}^{b} (a_h - n_c i_{g,1})^2$$

$$\leq \frac{8}{n^2(n-1)^2 b(b-1)} \sum_{g=1}^{b-1}\sum_{h=g+1}^{b} (a_h^2 - 2n_c g a_h + n_c^2 g^2)$$

$$\leq \frac{8}{n^2(n-1)^2 b(b-1)} \sum_{g=1}^{b-1}\left(\sum_{h=g+1}^{b}\left(\sum_{q=1}^{n_c}(n_c(h-1)+q)\right)^2 - 2n_c g \sum_{\tau=1}^{(b-g)n_c}(g+\tau) + n_c^2 g^2(b-g)\right) \quad (A.29)$$

$$= \frac{n_c^2\left((6n_c^2 - 2n_c - 2)b^2 + (2n_c^2 + 6n_c - 6)b + (-n_c^2 + 2n_c - 1)\right)}{3n^2(n-1)^2}$$

$$= \frac{(6b^2 + 2b - 1)n^2 + (-2b^2 + 6b + 2)bn - (2b^2 + 6b + 1)b^2}{3b^4(n-1)^2}$$

Since the distance is normalized to [0, 1], the upper bound of the variance of inter-neighbor distance is
$$VarD \leq \min\{1, \max ED^2 - \min E^2 D\}$$

$$\leq \min\left\{1, \frac{(6b^2 + 2b - 1)n^2 + (-2b^2 + 6b + 2)bn - (2b^2 + 6b + 1)b^2}{3b^4(n-1)^2} - \left(\frac{3n + (2b-1)b}{3b^2(n-1)}\right)^2\right\}.$$

$$= \min\left\{1, \frac{6(3b^2 + b - 2)n^2 - 6(b-2)(b+1)bn - 2(5b+2)(b+1)b^2}{9b^4(n-1)^2}\right\}$$

(A.30)

*Variance's lower bound*: The lower bound of the second moment of inter-neighbor distance is

$$ED^2 = \frac{8}{n^2(n-1)^2 b(b-1)} \cdot \frac{1}{\sum_{g=1}^{b-1}\sum_{h=g+1}^{b} 1^2} \cdot \sum_{g=1}^{b-1}\sum_{h=g+1}^{b}(a_h - n_c i_{g,1})^2 \cdot \sum_{g=1}^{b-1}\sum_{h=g+1}^{b} 1^2$$

$$\geq \frac{8}{n^2(n-1)^2 b(b-1)} \cdot \frac{2}{b(b-1)} \cdot \sum_{g=1}^{b-1}\left(\sum_{h=g+1}^{b}\sum_{q=1}^{n_c}(i_{h,q} - i_{g,1})\right)^2$$

$$\geq \frac{16}{n^2(n-1)^2 b^2(b-1)^2} \sum_{g=1}^{b-1}\left(\sum_{\tau=1}^{(b-g)n_c} \tau\right)^2 \quad (A.31)$$

$$= \frac{4n_c^2}{n^2(n-1)^2 b^2(b-1)^2} \sum_{g=1}^{b-1}\left((b + b^2 n_c) - (2bn_c + 1)g + n_c g^2\right)^2$$

$$= \frac{2n_c^2\left((6b^3 - 9b^2 + b + 1)n_c^2 + 15(b-1)bn_c + (10b-5)\right)}{15n^2(n-1)^2 b(b-1)}$$

$$= \frac{2\left((6b^3 - 9b^2 + b + 1)n^2 + 15(b-1)b^2 n + 5(2b-1)b^2\right)}{15(n-1)^2 b^5(b-1)}$$

Lower bound for the variance of inter-neighbor distance is
$$VarD \geq \max\{0, \min ED^2 - \max E^2 D\}$$

$$\geq \max\left\{0, \frac{2\left((6b^3 - 9b^2 + b + 1)n^2 + 15(b-1)b^2 n + 5(2b-1)b^2\right)}{15(n-1)^2 b^5(b-1)} - \left(\frac{(4b+1)n + (-2b+1)b}{3b^2(n-1)}\right)^2\right\}.$$

$$= \max\left\{0, \frac{(-80b^4 - 4b^3 - 89b^2 + b + 6)n^2 + 5(-16b^2 + 4b + 5)(b-1)b^2 n + 5(4b^4 - 8b^3 + 5b^2 + b - 1)b^2}{45(n-1)^2 b^5(b-1)}\right\}$$

(A.32)

# Appendix B: Processing Time for Operations on Eligible Machines

**Table 10** Processing time for operations on eligible machines

| Operation | $M_{1,1}$, $M_{1,2}$ (Centrifuge) | $M_{1,3}$ (Decapper) | $M_{1,4}$, $M_{1,5}$ (Biochemistry analyzer) | $M_{1,6}$, $M_{1,7}$ (Immunoassay analyzer) | $M_{1,8}$ (Validation /reporting) | $M_{2,1}$, $M_{2,2}$ (Centrifuge) | $M_{2,3}$ (Decapper) | $M_{2,4}$, $M_{2,5}$ (Biochemistry analyzer) | $M_{2,6}$, $M_{2,7}$ (Immunoassay analyzer) | $M_{2,8}$ (Validation /reporting) |
|---|---|---|---|---|---|---|---|---|---|---|
| $O_{1,1}$ | 545 | - | - | - | - | 332 | - | - | - | - |
| $O_{1,2}$ | - | 2 | - | - | - | - | 5 | - | - | - |
| $O_{1,3}$ | - | - | 545 | - | - | - | - | 530 | - | - |
| $O_{1,4}$ | - | - | - | - | 5 | - | - | - | - | 5 |
| $O_{2,1}$ | 593 | - | - | - | - | 356 | - | - | - | - |
| $O_{2,2}$ | - | 2 | - | - | - | - | 5 | - | - | - |
| $O_{2,3}$ | - | - | 593 | - | - | - | - | 695 | - | - |
| $O_{2,4}$ | - | - | - | - | 5 | - | - | - | - | 5 |
| $O_{3,1}$ | 530 | - | - | - | - | 325 | - | - | - | - |
| $O_{3,2}$ | - | 2 | - | - | - | - | 4 | - | - | - |
| $O_{3,3}$ | - | - | 530 | - | - | - | - | 475 | - | - |
| $O_{3,4}$ | - | - | - | - | 4 | - | - | - | - | 4 |
| $O_{4,1}$ | 597 | - | - | - | - | 358 | - | - | - | - |
| $O_{4,2}$ | - | 2 | - | - | - | - | 5 | - | - | - |
| $O_{4,3}$ | - | - | - | 1787 | - | - | - | - | 2669 | - |
| $O_{4,4}$ | - | - | - | - | 5 | - | - | - | - | 5 |
| $O_{5,1}$ | 516 | - | - | - | - | 318 | - | - | - | - |
| $O_{5,2}$ | - | 2 | - | - | - | - | 4 | - | - | - |
| $O_{5,3}$ | - | - | - | 1297 | - | - | - | - | 1442 | - |
| $O_{5,4}$ | - | - | - | - | 4 | - | - | - | - | 4 |
| $O_{6,1}$ | 564 | - | - | - | - | 342 | - | - | - | - |
| $O_{6,2}$ | - | 2 | - | - | - | - | 5 | - | - | - |
| $O_{6,3}$ | - | - | - | 1584 | - | - | - | - | 2161 | - |
| $O_{6,4}$ | - | - | - | - | 5 | - | - | - | - | 5 |

## Appendix C: Values for Decision Variables corresponding to a Feasible Schedule

**Table 11** Sets of decision variables taking value of 1. Variables not listed in this table take values of 0

| Job shop line #1 | $\{X_4^1, X_6^1\}$ |
| --- | --- |
| | $\{Y_{4,1,1}^{1,2}, Y_{4,2,1}^{1,3}, Y_{4,3,1}^{1,6}, Y_{4,4,1}^{1,8}, Y_{6,1,1}^{1,2}, Y_{6,2,2}^{1,3}, Y_{6,3,1}^{1,6}, Y_{6,4,2}^{1,8}\}$ |
| | $\{Z_1^{1,2,1}, Z_1^{1,3,2}, Z_2^{1,3,1}, Z_1^{1,6,1}, Z_1^{1,8,2}, Z_2^{1,8,1}\}$ |
| Job shop line #2 | $\{X_1^2, X_2^2, X_3^2, X_5^2\}$ |
| | $\{Y_{1,1,1}^{2,1}, Y_{1,2,1}^{2,3}, Y_{1,3,1}^{2,4}, Y_{1,4,1}^{2,8}, Y_{2,1,2}^{2,2}, Y_{2,2,2}^{2,3}, Y_{2,3,1}^{2,5}, Y_{2,4,2}^{2,8}, Y_{3,1,1}^{2,1}, Y_{3,2,3}^{2,3}, Y_{3,3,1}^{2,4}, Y_{3,4,3}^{2,8}, Y_{5,1,2}^{2,2}, Y_{5,2,4}^{2,3}, Y_{5,3,1}^{2,6}, Y_{5,4,4}^{2,8}\}$ |
| | $\{Z_1^{2,1,1}, Z_1^{2,2,1}, Z_1^{2,3,2}, Z_2^{2,3,4}, Z_3^{2,3,1}, Z_4^{2,3,3}, Z_1^{2,4,1}, Z_1^{2,5,1}, Z_1^{2,6,1}, Z_1^{2,8,2}, Z_2^{2,8,3}, Z_3^{2,8,1}, Z_4^{2,8,4}\}$ |

## Appendix D: Job Precedence Rule-based Distance Metric

To measure distance between solutions, Job Precedence Rule based distance measure, i.e., JPR distance (Liu et al., 2013) was adopted. JPR distance is the number of job pairs with identical elements but different precedence. Denote $\pi_1 = [3, 1, 6, 4, 5, 2]$ and $\pi_2 = [3, 4, 6, 1, 5, 2]$ as two permutations, and $D(\pi_1, \pi_2)$ as JPR distance between $\pi_1$ and $\pi_2$. For instance, we extract all pairs of jobs with precedence from $\pi_1$ and $\pi_2$. For $\pi_1$, we extracted (3,1), (3,6), (3,4), (3,5), (3,2), (1,6), (1,4), (1,5), (1,2), (6,4), (6,5), (6,2), (4,5), (4,2), and (5,2). For $\pi_2$, we extracted (3,4), (3,6), (3,1), (3,5), (3,2), (4,6), (4,1), (4,5), (4,2), (6,1), (6,5), (6,2), (1,5), (1,2), and (5,2). For the job pairs of 1 and 4, 1 precedes 4 in $\pi_1$, while 1 succeeds 4 in $\pi_2$, which adds a distance of 1. Similarly another two job pairs (i.e., 4 and 6), 4 and 6 show different precedence in $\pi_1$ and $\pi_2$. Thus, their JPR distance is 3. We normalized the JPR as $DIST(\pi_1, \pi_2) = D(\pi_1, \pi_2) / D_{\max}$, where $D_{\max} = n(n-1)/2$. Finally, the normalized JPR distance $DIST(\pi_1, \pi_2) = 3/(6 \times 5/2) = 0.2$.

## Appendix E: Concepts in Local Optima Network

**Table 12** Concepts in local optima network

| Concepts | Explanation |
| --- | --- |
| Node | Represents local optimum. |
| Escape edge | Escape edge exists between two nodes if a better local optimum can be obtained by performing search from the starting point of a local optimum. Escape edge is directed from the starting node to the ending node. |
| Edge weight | For escape edge, it is the frequency at which the directed transition between connected nodes occurs. |
| Edge width | Proportional to edge weight. |
| Node size | Proportional to weighted incoming degree, reflecting the node attraction. |
| Funnel structure | A group of local optima that converge to a single solution. |
| Big-valley structure (globally convex structure) | Local optima are clustered around a global optimum. It is a special case of funnel structure with only one funnel. |
| Plateau | A higher-level virtual node compressing a group of local optima with same objective values belonging to a connected component according to the escape edge. It can be further divided into attractor plateaus and sink plateaus. |
| Attractor plateau | A plateau that can be improved. |
| Sink plateau (Funnel bottom) | A plateau with the best objective value in the funnel. |
| Plateau average size | Average number of local optima within each plateau. |

## Appendix F: Recording Nodes and Edges for Local Optima Network

Two basic elements for a visualized local optima network (LON) are nodes and edges. Nodes and edges found during the search on an instance by an optimization algorithm were recorded. For generating an LON, each optimization algorithm is independently repeated $N$ times on each instance. The algorithm is terminated if the best-so-far solution has been unimproved after $M$ iterations. We use node set $L$ to store all newly-found improved local optima, and edge set $E$ to record all transitions starting from a local optimum and ending at another optimum. If a local optimum/edge occurs multiple times, duplicates were deleted; only one was kept. The incoming degree for a node (node size in LON) was recorded, as well as frequency (edge width in LON) at which the directed transition between connected nodes occurs. Two parameters were set: for toy instances, the settings in Ochoa and Veerapen (2018), meaning that the number of independent runs $N=1000$ and maximum number of iterations $M=10,000$; for medium- and large-size instances, $N=50$ and $M=1000$ were used to save computational cost.

**Author contributions** Keyao Wang developed the method, performed the experiments, and contributed to paper writing. Bo Liu raised the research funding, initiated the idea, supervised the first author, and wrote the paper.

**Funding** Bo Liu was funded by Frontier Science Key Research Program, Chinese Academy of Sciences under Grant No. QYZDB-SSW-SYS020.

**Data availability** All benchmark instances, datasets and a guide on how to use it are made publicly available at https://github.com/kywang95/HDFJSP.

**Code availability** We agree upon publishing codes in case of acceptance.

## Declarations

**Conflict of interest** The authors declare that they have no conflict of interest.

**Ethics approval** Not Applicable.

**Consent to participate** Not Applicable.

**Consent to publication** Not Applicable.

# References

Aarts, R. J., Lindsey, J. S., Corkan, L. A., & Smith, S. F. (1995). Flexible Protocols Improve Parallel Experimentation Throughput. *Clinical Chemistry, 41*(7), 1004-1010.
Ball, R. C., Branke, J., & Meisel, S. (2018). Optimal Sampling for Simulated Annealing Under Noise. *INFORMS Journal on Computing, 30*(1), 200-215.
Barbulescu, L., Howe, A. E., Whitley, L. D., & Roberts, M. (2006). Understanding algorithm performance on an oversubscribed scheduling application. *Journal of Artificial Intelligence Research, 27*, 577-615.
Barrett, A., & Weld, D. S. (1994). Partial-order planning: Evaluating possible efficiency gains. *Artificial Intelligence, 67*(1), 71-112.
Behnamian, J., & Ghomi, S. M. T. F. (2016). A survey of multi-factory scheduling. *Journal of Intelligent Manufacturing, 27*(1), 231-249.
Bengio, Y., Lodi, A., & Prouvost, A. (2021). Machine learning for combinatorial optimization: A methodological tour d'horizon. *European Journal of Operational Research, 290*(2), 405-421.
Boyd, J. C., & Savory, J. (2001). Genetic algorithm for scheduling of laboratory personnel. *Clinical Chemistry, 47*(1), 118-123.
Chang, H. C., & Liu, T. K. (2017). Optimisation of distributed manufacturing flexible job shop scheduling by using hybrid genetic algorithms. *Journal of Intelligent Manufacturing, 28*(8), 1973-1986.
Daolio, F., Liefooghe, A., Verel, S., Aguirre, H., & Tanaka, K. (2017). Problem Features versus Algorithm Performance on Rugged Multiobjective Combinatorial Fitness Landscapes. *Evolutionary Computation, 25*(4), 555-585.
Daolio, F., Verel, S., Ochoa, G., & Tomassini, M. (2014). Local Optima Networks of the Permutation Flow-Shop Problem. In P. Legrand, M.-M. Corsini, J.-K. Hao, N. Monmarché, E. Lutton & M. Schoenauer (Eds.), *Artificial Evolution: 11th International Conference, Evolution Artificielle, EA 2013, Bordeaux, France, October 21-23, 2013. Revised Selected Papers* (pp. 41-52). Cham: Springer International Publishing.
Davis, J. (1995). Scheduling in the Clinical Laboratory. *Clinical Chemistry, 41*(7), 961-962.
De Giovanni, L., & Pezzella, F. (2010). An Improved Genetic Algorithm for the Distributed and Flexible Job-shop Scheduling problem. *European Journal of Operational Research, 200*(2), 395-408.
Drake, J. H., Kheiri, A., Özcan, E., & Burke, E. K. (2020). Recent advances in selection hyper-heuristics. *European Journal of Operational Research, 285*(2), 405-428.
Du, Y., Li, J. Q., Luo, C., & Meng, L. L. (2021). A hybrid estimation of distribution algorithm for distributed flexible job shop scheduling with crane transportations. *Swarm and Evolutionary Computation, 62*.
Elena, B. M., Beraldi, P., & Conforti, D. (2006). Improving the efficiency of a clinical laboratory: a mathematical approach. *IFAC Proceedings Volumes, 39*(3), 659-664.
Esposito, L. (2015). Hospital Labs: Behind the Scenes. <https://health.usnews.com/health-news/patient-advice/articles/2015/01/30/hospital-labs-behind-the-scenes> [accessed 03/09/2018].
Fernandez-Viagas, V., & Framinan, J. M. (2015). A bounded-search iterated greedy algorithm for the distributed permutation flowshop scheduling problem. *International Journal of Production Research, 53*(4), 1111-1123.
Franzin, A., & Stützle, T. (2023). A landscape-based analysis of fixed temperature and simulated annealing. *European Journal of Operational Research, 304*(2), 395-410.
Frost, D., & Dechter, R. (1994). In Search of the Best Constraint Satisfaction Search. *Proceedings of the Twelfth National Conference on Artificial Intelligence, Vols 1 and 2*, 301-306.
Glover, F. (1994). Tabu Search for Nonlinear and Parametric Optimization (with Links to Genetic Algorithms). *Discrete Applied Mathematics, 49*(1-3), 231-255.
Glover, F. (1998). A template for scatter search and path relinking. In J. Hao, E. Lutton, E. Ronald, M. Schoenauer & S. D. (Eds.), *Artificial Evolution* (Vol. 1363, pp. 1-51). Berlin, Heidelberg: Springer.
Goldberg, D. E. (1989). *Genetic Algorithms in Search, Optimization and Machine Learning*: Addison-Wesley Longman Publishing Co., Inc.
Gratch, J., & Chien, S. (1996). Adaptive problem-solving for large-scale scheduling problems: A case study. *Journal of Artificial Intelligence Research, 4*, 365-396.
Gratch, J., & DeJong, G. (1996). A statistical approach to adaptive problem solving. *Artificial Intelligence, 88*(1), 101-142.
Hatami, S., Ruiz, R., & Andres-Romano, C. (2013). The Distributed Assembly Permutation Flowshop Scheduling Problem. *International Journal of Production Research, 51*(17), 5292-5308.
Holland, L. L., Smith, L. L., & Blick, K. E. (2006). Total laboratory automation can help eliminate the laboratory as a factor in emergency department length of stay. *American Journal of Clinical Pathology, 125*(5), 765-770.
Hooker, J. N. (2007). Planning and scheduling by logic-based benders decomposition. *Operations Research, 55*(3), 588-602.
Hutter, F., Xu, L., Hoos, H. H., & Leyton-Brown, K. (2014). Algorithm runtime prediction: Methods & evaluation. *Artificial Intelligence, 206*, 79-111.
Jones, T., & Forrest, S. (1995). *Fitness Distance Correlation as a Measure of Problem Difficulty for Genetic Algorithms*. Paper presented at the Proceedings of the Sixth International Conference on Genetic Algorithms.
Kambhampati, S., Knoblock, C. A., & Yang, Q. (1995). Planning as Refinement Search - a Unified Framework for Evaluating Design Tradeoffs in Partial-Order Planning. *Artificial Intelligence, 76*(1-2), 167-238.
Kerschke, P., Hoos, H. H., Neumann, F., & Trautmann, H. (2019). Automated Algorithm Selection: Survey and Perspectives. *Evolutionary Computation, 27*(1), 3-45.
Keskinocak, P., & Savva, N. (2020). A Review of the Healthcare-Management (Modeling) Literature Published in Manufacturing & Service Operations Management. *Manufacturing & Service Operations Management, 22*(1), 59-72.
Kim, J. W., Choi, B. J., Noh, K. H., Choi, H. R., Koo, J. C., Ryew, S. M., et al. (2006). *Automatic scheduling algorithm for personalized clinical test*. Paper presented at the 2006 SICE-ICASE International Joint Conference.
Kirkpatrick, S., Gelatt, C. D., & Vecchi, M. P. (1983). Optimization by Simulated Annealing. *Science, 220*(4598), 671-680.
Kool, W., van Hoof, H., & Welling, M. (2019). *Attention, Learn to Solve Routing Problems!* Paper presented at the ICLR.
Liu, B., Wang, L., & Jin, Y. H. (2007). An effective PSO-based memetic algorithm for flow shop scheduling. *IEEE Transactions on Systems Man and Cybernetics Part B-Cybernetics, 37*(1), 18-27.
Liu, B., Wang, L., Liu, Y., & Wang, S. Y. (2011). A unified framework for population-based metaheuristics. *Annals of Operations Research, 186*(1), 231-262.
Liu, B., Xu, J., Qian, B., Wang, J., & Chu, Y. (2013, 16-19 April 2013). *Probabilistic memetic algorithm for flowshop scheduling*. Paper presented at the 2013 IEEE Workshop on Memetic Computing (MC).
Lofberg, J. (2005). *A toolbox for modeling and optimization in MATLAB*. Paper presented at the IEEE International Symposium on Computer Aided Control Systems Design.
Lohmer, J., & Lasch, R. (2020). Production planning and scheduling in multi-factory production networks: a systematic literature review. *International Journal of Production Research*, 1-27.
Lu, P. H., Wu, M. C., Tan, H., Peng, Y. H., & Chen, C. F. (2018). A genetic algorithm embedded with a concise chromosome representation for distributed and flexible job-shop scheduling problems. *Journal of Intelligent Manufacturing, 29*(1), 19-34.
Malan, K. M., & Engelbrecht, A. P. (2013). A survey of techniques for characterising fitness landscapes and some possible ways forward. *Information Sciences, 241*, 148-163.
Mendez, C. A., Cerda, J., Grossmann, I. E., Harjunkoski, I., & Fahl, M. (2006). State-of-the-art review of optimization methods for short-term scheduling of batch processes. *Computers & Chemical Engineering, 30*(6-7), 913-946.
Meng, L. L., Zhang, C. Y., Ren, Y. P., Zhang, B., & Lv, C. (2020). Mixed-integer linear programming and constraint programming formulations for solving distributed flexible job shop scheduling problem. *Computers & Industrial Engineering, 142*.
Merz, P., & Katayama, K. (2004). Memetic algorithms for the unconstrained binary quadratic programming problem. *Biosystems, 78*(1-3), 99-118.
Metropolis, N., Rosenbluth, A. W., Rosenbluth, M. N., Teller, A. H., & Teller, E. (1953). Equation of State Calculations by Fast Computing Machines. *The Journal of Chemical Physics, 21*(6), 1087-1092.
Minton, S. (1988). *Learning search control knowledge: an explanation-based approach*. Boston: Kluwer Academic.
Naderi, B., & Azab, A. (2014). Modeling and heuristics for scheduling of distributed job shops. *Expert Systems with Applications, 41*(17), 7754-7763.
Naderi, B., & Ruiz, R. (2014). A scatter search algorithm for the distributed permutation flowshop scheduling problem. *European Journal of Operational Research, 239*(2), 323-334.
Nawaz, M., Enscore, E. E., & Ham, I. (1983). A Heuristic Algorithm for the M-Machine, N-Job Flowshop Sequencing Problem. *Omega-International Journal of Management Science, 11*(1), 91-95.
Nguyen, Q. H., Ong, Y. S., & Lim, M. H. (2009). A Probabilistic Memetic Framework. *IEEE Transactions on Evolutionary Computation, 13*(3), 604-623.
Nowicki, E., & Smutnicki, C. (2006). Some aspects of scatter search in the flow-shop problem. *European Journal of Operational Research, 169*(2), 654-666.
Ochoa, G., & Veerapen, N. (2018). Mapping the global structure of TSP fitness landscapes. *Journal of Heuristics, 24*(3), 265-294.
Ong, Y. S., & Keane, A. J. (2004). Meta-Lamarckian learning in memetic algorithms. *IEEE Transactions on Evolutionary Computation, 8*(2), 99-110.
Pezzella, F., Morganti, G., & Ciaschetti, G. (2008). A genetic algorithm for the Flexible Job-shop Scheduling Problem. *Computers & Operations Research, 35*(10), 3202-3212.
Pinedo, M. L. (2012). *Scheduling: theory, algorithms, and systems* (4 ed.). Boston, MA: Springer.
Reeves, C. R. (1999). Landscapes, operators and heuristic search. *Annals of Operations Research, 86*, 473-490.
Reidys, C. M., & Stadler, P. F. (2002). Combinatorial landscapes. *SIAM Review, 44*(1), 3-54.
Ribas, I., Companys, R., & Tort-Martorell, X. (2010). Comparing three-step heuristics for the permutation flow shop problem. *Computers & Operations Research, 37*(12), 2062-2070.
Rice, J. R. (1976). The Algorithm Selection Problem. In M. Rubinoff & M. C. Yovits (Eds.), *Advances in Computers* (Vol. 15, pp. 65-118): Elsevier.
Roshanaei, V., Luong, C., Aleman, D. M., & Urbach, D. R. (2017). Collaborative Operating Room Planning and Scheduling. *INFORMS Journal on Computing, 29*(3), 558-580.
Rudin, C., Chen, C. F., Chen, Z., Huang, H. Y., Semenova, L., & Zhong, C. D. (2022). Interpretable machine learning: Fundamental principles and 10 grand challenges. *Statistics Surveys, 16*, 1-85.
Schmidt, G. (1998). Case-based reasoning for production scheduling. *International Journal of Production Economics, 56-57*, 537-546.
Singer, A. J., Ardise, J., Gulla, J., & Cangro, J. (2005). Point-of-care testing reduces length of stay in emergency department chest pain patients. *Annals of Emergency Medicine, 45*(6), 587-591.
Streeter, M. J., & Smith, S. F. (2006). How the landscape of random job shop scheduling instances depends on the ratio of jobs to machines. *Journal of Artificial Intelligence Research, 26*, 247-287.
Toptal, A., & Sabuncuoglu, I. (2010). Distributed scheduling: a review of concepts and applications. *International Journal of Production Research, 48*(18), 5235-5262.
Van Merode, G. G., Oosten, M., Vrieze, O. J., Derks, J., & Hasman, A. (1998). Optimisation of the structure of the clinical laboratory. *European Journal of Operational Research, 105*(2), 308-316.
Wang, J. J., & Wang, L. (2020). A Knowledge-Based Cooperative Algorithm for Energy-Efficient Scheduling of Distributed Flow-Shop. *IEEE Transactions on Systems Man Cybernetics-Systems, 50*(5), 1805-1819.
Wang, K., & Liu, B. (2022). HDFJSP. <https://github.com/kywang95/HDFJSP>. *GitHub repository*.
Watson, J. P., Barbulescu, L., Whitley, L. D., & Howe, A. E. (2002). Contrasting structured and random permutation flow-shop scheduling problems: Search-space topology and algorithm performance. *Informs Journal on Computing, 14*(2), 98-123.
Watson, J. P., Beck, J. C., Howe, A. E., & Whitley, L. D. (2003). Problem difficulty for tabu search in job-shop scheduling. *Artificial Intelligence, 143*(2), 189-217.
Watson, J. P., Whitley, L. D., & Howe, A. E. (2005). Linking search space structure, run-time dynamics, and problem difficulty: A step toward demystifying tabu search. *Journal of Artificial Intelligence Research, 24*(1), 221-261.
Weinberger, E. (1990). Correlated and Uncorrelated Fitness Landscapes and How to Tell the Difference. *Biological Cybernetics, 63*(5), 325-336.
Wolpert, D. H., & Macready, W. G. (1997). No free lunch theorems for optimization. *IEEE Transactions on Evolutionary Computation, 1*(1), 67-82.
Yang, Y., Li, P., Wang, S., Liu, B., & Luo, Y. (2017). *Scatter search for distributed assembly flowshop scheduling to minimize total tardiness*. Paper presented at the 2017 IEEE Congress on Evolutionary Computation (CEC).
Zhang, W. (2004). Configuration landscape analysis and backbone guided local search.: Part I: Satisfiability and maximum satisfiability. *Artificial Intelligence, 158*(1), 1-26.
Ziaee, M. (2014). A heuristic algorithm for the distributed and flexible job-shop scheduling problem. *Journal of Supercomputing, 67*(1), 69-83.